\documentclass{article}
\usepackage[margin=1in]{geometry}
\usepackage[utf8]{inputenc}
\usepackage{amsthm}
\pdfoutput=1 
\usepackage{amsmath}
\usepackage{graphicx}
\usepackage[space]{grffile}
\usepackage{hyperref}
\usepackage{caption}
\usepackage{subcaption}
\usepackage{bm}
\usepackage{import}
\usepackage{kantlipsum}
\usepackage{enumitem}
\usepackage[utf8]{inputenc}
\usepackage[english]{babel}
\newtheorem{theorem}{Theorem}[section]
\usepackage{booktabs}
\usepackage{multirow}
\usepackage{placeins}
\usepackage{cite}
\usepackage[T1]{fontenc}

\usepackage{xcolor}
\usepackage{ulem}
\usepackage{amsfonts}
\usepackage{todonotes}
\usepackage{bm}
\newcommand{\R}{\mathbb{R}}
\newcommand{\T}{\Omega}
\newcommand{\Pu}{\mathcal{P}_u}
\newcommand{\nn}{\mathcal{NN}}
\newcommand{\x}{\mathbf{x}}
\DeclareMathOperator*{\argmax}{arg\,max}
\DeclareMathOperator*{\argmin}{arg\,min}
\newcommand{\NNu}{\NN_{u}}

\title{Neural Network Based Variational Methods for Solving Quadratic Porous Medium Equations in High Dimensions \footnote{This work is supported in part by National Science Foundation via grant DMS-2012286 and by Department of Energy via grant DE-SC0019449. }}
\author{Jianfeng Lu, Min Wang}

\begin{document}

\maketitle
\begin{abstract}
In this paper, we propose and study neural network based methods for solutions of high-dimensional quadratic porous medium equation (QPME). Three variational formulations of this nonlinear PDE are presented: a strong formulation and two weak formulations. For the strong formulation, the solution is directly parameterized with a neural network and optimized by minimizing the PDE residual. It can be proved that the convergence of the optimization problem guarantees the convergence of the approximate solution in the $L^1$ sense.
The weak formulations are derived following \cite{brenier2020examples} which characterizes the very weak solutions of QPME. 
Specifically speaking, the solutions are represented with intermediate functions who are parameterized with neural networks and are trained to optimize the weak formulations. 
Extensive numerical tests are further carried out to investigate the pros and cons of each formulation in low and high dimensions. This is an initial exploration made along the line of solving high-dimensional nonlinear PDEs with neural network based methods, which we hope can provide some useful experience for future investigations.
\end{abstract}

\section{Introduction}

Solving high-dimensional PDEs is a long-standing challenge in scientific computing. Standard mesh-based methods, such as Finite Element Method (FEM), Finite Difference Method (FDM) would suffer from the curse of dimensionality, i.e., in order to sustain the accuracy of the approximate solution to high-dimensional PDEs, an approximation space with exponentially large size must be used. The number of degrees of freedom associated to the approximation space is often proportional to the number of elements in the mesh which usually scales exponentially in the dimension to achieve an suitable discretization of the domain. Therefore, for high-dimensional problems, the mesh-based methods are impractical. Alternatively, semilinear parabolic PDEs can also be solved point-wisely based on stochastic representation of the solutions using Monte Carlo algorithm \cite{henry2012counterparty,warin2018nesting,henry2019branching,warin2017variations}, but such approaches only apply to specific type of PDEs.

To circumvent the challenges for solving general high dimensional (nonlinear) PDEs, many attempts have been made. One natural idea is to restrict to a solution ansatz. For example, using the tensor train (TT) format to approximate the solutions of high-dimensional PDEs \cite{richter2021solving,dektor2021rank,dektor2020dynamically,dolgov2021tensor,eigel2019non,boelens2018parallel}. While such methods are quite successful if the solution can be well represented by the tensor train, the representability is not guaranteed. Another natural and promising candidate for PDE solution ansatz is the artificial neural networks. Thanks to the rich expressiveness of the neural networks to parametrize high dimensional functions \cite{barron1993universal}. Theoretical results are also available to justify the approximability of PDE solutions by neural networks without curse of dimensionality, e.g.,    \cite{jentzen2018proof,lu2021priori}.
Many recent works have been devoted for various approaches of using neural networks to solve high dimensional PDEs. Typically, such methods first identify functional optimization problems corresponding to the PDEs
$$\underbrace{u}_{\text{PDE solution}}   =  \underbrace{\argmin_f\mathcal{C}(f)}_{\text{PDE-inspired optimization problem}}.$$
Then one could take neural network  as an ansatz of the minimizer $u$ and minimize the parameters using stochastic gradient type approaches. 

One well known method of this kind is the physics informed neural network (PINN) \cite{raissi2019physics}, which takes the loss function to be directly the PDE residual.
One drawback of PINN is that, to compute the loss function, the derivatives or high-order derivatives of the neural network need to be computed throughout the training process. However, the generality and simplicity of this framework still make it an easy choice when it comes to the high-dimensional PDEs. 
Other neural network based methods to solve PDEs include the Deep Ritz Methods \cite{yu2018deep},  the Deep BSDE method  \cite{han2018solving,han2017deep} for semilinear parabolic PDEs, utilizing different formulations to turn the high dimensional PDE problems into an optimization problem on the neural network parameters. 
More recent work \cite{zang2020weak} adopts weak variational formulation for linear equations to solve such PDEs using neural networks. While it is generically unclear how to extend such techniques to nonlinear ones. 

In this paper, we make an attempt along the direction of utilizing weak formulation for solving nonlinear PDEs in high dimension. In particular, we consider high-dimensional quadratic porous medium equation (QPME) (Section \ref{QPME_intro}). Several variational formulations of the QPME is proposed to solve such PDE (Section \ref{sec:formulation}) which both allow solutions in a very weak sense and can be sought with deep learning techniques.
In addition, these formulations are further compared with the PINN formulation. In Section \ref{sec:nn}, more detailed treatment of the neural network for solving the optimization problem is presented. Numerical results are then provided to verify and compare the effectiveness of the proposed methods in Section \ref{sec:nuemrical} and Section \ref{sec:nuemrical_waiting}.

\section{Preliminaries}\label{QPME_intro}
We consider the porous medium equation (PME)
$$\partial_t u =  \frac{1}{m}\Delta u^m,\quad (t,\x) \in Q.$$
PME is a degenerate parabolic equation.  It is only parabolic where $u>0$. 
When $m$ is taken to be $2$, the quadratic porous medium equation (QPME or Boussinesq's equation) reads 
\begin{equation}\label{QPME}
\begin{split}
    \partial_t u = \frac{1}{2}\Delta u^2 = \text{div}\left(u \nabla u \right)  = \frac{1}{2} u\Delta u + \frac{1}{2}|\nabla u|^2 ,\quad (t,\x) \in Q,
\end{split}
\end{equation}
where $Q := [0,T]\times \T$ and $\T$ stands for a bounded domain in $\R^d$. $\Delta= \Delta_x$ represents the Laplace operator acting on the space variables.
The equation is mainly used to describe process involving fluid flow, heat transfer or diffusion. Particularly, it can be used to study the flow of an isentropic gas through a porous medium \cite{muskat1938flow,leibenzon1930motion}. In this case, $u(t,\x)\in \R$ is a scalar function denoting the density of the gas ($u^{2}$ is roughly the pressure, and $u\nabla u$ stands for the flux). Physical constraints may apply such that $u> 0$. Power $2$ here relates to the thermal dynamics character expansion in terms of the pressure of the gas of interest (linear dependency).

A main feature of PME is the ``finite propagation'' speed, in contrast with the infinite propagation of usual diffusion equations. Essentially, the free boundary that separates the regions where the solution is positive (i.e. where ``there is gas", according to the standard interpretation of $u$ as a gas density), from the ``empty region" where $u = 0$ moves as time passes.
$$
\Gamma  = \partial \Pu \cap Q
$$
where $\Pu := \{(t,\x)\in Q\ |\ u(t,\x)>0 \}$ denotes the set where $u$ is positive.
$\Gamma$ is also sometimes referred as the moving boundary, propagation fronts, or the interface. While a rather comprehensive theoretical analysis of this PDEs is provided in \cite{vazquez2007porous}, exact solutions to general initial/boundary conditions can usually not be obtained. Numerical schemes thus must be applied to obtain approximate solutions. Most previous studies of PME from a numerical aspect put their focus on dealing with the moving free boundaries of the solutions. The adaptive moving mesh schemes were proposed and coupled with mesh-based methods such as finite element method (FEM) to obtain accurate yet efficient numerical solutions to PME \cite{ngo2017study}. However, such methods can not be used to solve high-dimensional QPME due to the curse of dimensionality. The only exception as far as we know is \cite{shukla1996use}, in which the supervised learning was conducted to learn the correspondence between the physical parameters in the PDE and the one-dimensional solutions at certain evaluation $x$, the learning of the global solution is not considered by the authors. 

While PME is mainly used in modellings for low physical dimensions, $d = 2$ or $3$. In this work, we use it as a prototypical degenerate nonlinear equation in high dimensions to test numerical PDE solvers based on the neural networks. The high dimensional diffusion might be used for certain machine learning tasks such as analysis of high dimensional point cloud data, which we will leave for future investigations. 


\section{Variational formulations of QPME}\label{sec:formulation}
 Since the mesh-based algorithms suffer from the curse of dimensionality, we therefore turn to neural network based  techniques for solutions to high dimensional PDEs. 
In particular, we first convert the  initial/boundary value problem (IBVP) of QPME into a variational formulation  and then take the neural network as an ansatz of the solution. The objective function is then taken as the loss function and the extrema will be obtained by optimizing the loss function with stochastic gradient descent (SGD) or its variants.

In this section, we specifically focus on the first step of this procedure, i.e., the IBVP and the variational reformulation. 

\subsection{Initial / boundary value problem}

Consider the QPME on a hyperrectangle
$$\partial_t u= \frac{1}{2} \Delta u^2,  \quad (t,\x )\in Q$$
where $Q =[0,T]\times \T,$ and $ \T = [-a_i,a_i]^{d}$. We consider the QPME with 
the homogeneous Dirichlet boundary condition 
\begin{equation}\label{BC}
    \textbf{Dirichlet B.C.}\quad  u(t,\x)|_{\Sigma_T} = 0
\end{equation}
where $\Sigma_T: = [0,T] \times\partial \T$. We also impose the initial condition to the PDE as 
\begin{equation}\label{IC}
   \textbf{I.C.}\quad  u(0,\x) = u_0(\x)\quad \x\in \T.
\end{equation}

\subsection{Strong formulation}\label{sec:PINN}
One immediate optimization formulation is to use the strong form of the PDE by minimizing the squared PDE residual
\begin{equation}\label{PINN}
    \mathcal{L}_{\text{PDE}} (u) = \int_Q \left( \partial _t u -\frac{1}{2} \Delta u^2\right)^2.
\end{equation}
If both the I.C.{} and  B.C.{} are strictly enforced as hard constraints, the optimization problem can then be formulated as 
\begin{equation}\label{PINNN_strong}
    \min_{u\in V_0} \mathcal{L}_{\text{PDE}} (u)
\end{equation}
where $V_0: = \{ f :  f|_{\Sigma_t} = 0 , f(0,\x) = u_0(\x)\}$.
Alternatively, both  I.C.{} and  B.C.{} can be treated as soft constraints  enforced by penalizations:  We may define
\begin{equation}\label{bc_weak}
    \mathcal{L}_{B}(u)  = \int_{\Sigma_T} u ^2
\end{equation}
for homogeneous Dirichlet boundary condition, and
\begin{equation}\label{weak_initial}
      \mathcal{L}_{I}(u) := \int_{\T} \left(u \left(0,\x\right) - u_{0}\left(\x\right)\right)^2 
\end{equation}
for the initial condition.
The optimization problem \eqref{PINNN_strong} can then be relaxed to 
\begin{equation}\label{PINN_full}
    \min_{u\in V} \mathcal{L}_{\text{PINN}}(u)
\end{equation}
for some function space $V$, where 
\begin{equation}   \mathcal{L}_{\text{PINN}}(u): = \kappa  \mathcal{L}_{\text{PDE}}(u) + \mu\mathcal{L}_{B}(u) +
    \nu\mathcal{L}_{I}(u)
\end{equation}
is weighted sum of the PDE residual, the error of boundary condition and the error for initial condition. $\kappa, \mu,\nu$ are weights for each term. 
We used the subscript PINN for the loss function, as such formulation was popularized by the PINN method \cite{raissi2019physics} in recent years, while the idea dates back to early days of using neural network ansatz for PDE solutions \cite{lagaris1998artificial}. 

So far, the PDE residual, mismatch in initial condition and boundary condition of $u$ are all measured in a $L^2$ sense. We could also define an analogous $L^1$ optimization problem \eqref{PINN_full} with:
\begin{equation}\label{PINN_L1}
\begin{split}
    &\mathcal{L}_{\text{PDE}} (u) = \int_Q \left|\partial _t u -\frac{1}{2} \Delta u^2\right|,
    \\
     & 
     \mathcal{L}_{B}(u)  = \int_{\Sigma_T} |u|,\\
     &\mathcal{L}_{I}(u;u_0) := \int_{\T} |u \left(0,\x\right) - u_{0}\left(\x\right)|.
\end{split}
\end{equation}
We refer the target function in $L^1$ by $\mathcal{L}_{\text{PINN}-L^1}$ and that in $L^2$ by $\mathcal{L}_{\text{PINN}-L^2}$.
While using $L^2$ to measure the PDE residual and I.C./B.C. mismatch is a standard practice in PINN,  the use of $L^1$ is inspired by the following stability analysis.


Assume $u$ is a strong solution to the homogeneous QPME \eqref{QPME} subject to the homogeneous Dirichlet boundary condition and the initial condition \eqref{IC}. Let $\hat{u}(t, \x)$ be another smooth function that satisfies the Dirichlet boundary condition $\hat{u}(t,\x)|_{\Sigma_t} = 0$. We define $\hat{f}(t,\x) := \partial_t \hat{u}(t,\x) - \frac{1}{2}\Delta \hat{u}^{2} (t,\x)$, so that $\hat{u}$ satisfies the QPME with $\hat{f}$ being the forcing term
\begin{equation*}
    \begin{split}
        &\partial_t \hat{u} = \frac{1}{2} \Delta \hat{u}^2 +\hat{f},\quad \forall (t,\x)\in Q\\
        &\hat{u}(0,\x;\theta)= \hat{u}_0(\x), \quad \forall \x \in \T,
    \end{split}
\end{equation*}
where $\hat{u}_0(\x) = \hat{u}(0, \x)$.

Then the following $L^1$ contraction holds  for any $t\in[0,T]$  \cite{vazquez2007porous}:
    \begin{equation}
        ||u(t)- \hat{u}(t)||_{1}\leq ||u_0 - \hat{u}_0||_1 + \int_{0}^t ||f(s) - \hat{f}(s)||_1\ ds.
    \end{equation}
Thus, 
\begin{equation}\label{convergence_L1}
    ||u(t)- \hat{u}(t)||_{1}\leq ||u_0 - \hat{u}_0||_1 + \int_{0}^t ||\partial_t \hat{u} - \frac{1}{2}\Delta \hat{u}^2 ||_1\ ds\leq \tilde{C} \mathcal{L}_{\text{PINN}-L^1} (\hat{u})
\end{equation}
noticing $f(t)$ in \eqref{QPME} is zero.
Therefore, as long as $\T$ is a bounded domain, by Cauchy–Schwarz inequality, we have
\begin{equation} \label{convergence_L2}
    ||u(t)- \hat{u}(t)||^2_{1}\leq C\big( ||u_0 - \hat{u}_0||^2_2+ \int_{0}^t ||\partial_t \hat{u} - \frac{1}{2}\Delta \hat{u}^2 ||^2_2\ ds \big)\leq \tilde{C} \mathcal{L}_{\text{PINN}-L^2} (\hat{u}).
\end{equation}
The estimates 
\eqref{convergence_L1} and \eqref{convergence_L2} show that if the QPME is solved by minimizing the loss $\mathcal{L}_{\text{PINN}}$, the convergence of the optimization guarantees the convergence to the exact solution to \eqref{QPME} in $L^1$ sense. In other words, if the problem is further solved with deep learning methods, the loss is essentially an direct real time  indicator of the approximation error during the entire training process. It can thus be naturally used to truncate the training when needed. Additionally, this result suggests that the intrinsic norm to measure the approximation error is $L^1$ norm.

\subsection{Relaxed concave optimization problem}
Besides the strong formulation, in this section, we derive and consider a series of optimization problems which can also be used to solve the QPME; they correspond to various weak formulation of the PDE. 

We start by considering the {very weak solutions} to the QPME \eqref{QPME}, i.e. $u \in L^{1}(Q)$ satisfying 
\begin{equation}
    \int_{Q } -2\partial_t\psi u -\Delta \psi u^2 +2u_0 \partial_t \psi =0\quad 
\end{equation}
for all test functions $\psi \in C^{2,1}(\bar{Q})$ which vanishes on $\Sigma_T$ and for $t= T$. Essentially, a very weak solution is an integrable distribution solution. Unlike strong solutions, no derivative of the solution is used in defining such solutions; so very weak solutions have much lower regularity requirements. 
We also remark that while we focus on very weak solutions in this paper, there are different ways of defining generalized solutions for QPME. A weak solution, for example, is defined to be a function $u$ such that $u^2\in L_{loc}^1(0,T;\ W_{loc}^{1,1})$ which satisfies 
$$\int_Q -2\partial_t\psi u +\nabla  (u^2)\cdot\nabla \psi + 2u_0\partial_t \psi = 0$$
It is clear that all weak solutions are very weak solutions by definition; weak solutions require higher regularity of the solutions.

The following theorem gives a characterization  for very weak solutions to QPME \cite{brenier2020examples}. 
\begin{theorem}[\cite{brenier2020examples}]\label{thm}
        Any \textbf{very weak solution} $u$ to QPME can be recovered as 
        \begin{equation}\label{u_phi_thm}
             u = \frac{\partial_t \phi^*}{1-\Delta \phi^*}
        \end{equation}
        where 
        \begin{equation}\label{phi_formulation}
          \phi^*  =  {\argmax_{\phi\in B} J(u_0)} = \argmax_{\phi\in B} \int_Q \frac{-(\partial_t \phi)^2}{ 1-\Delta \phi} + 2u_0 \partial_t \phi
        \end{equation}
        with $B:= \{\phi \ |\  \phi(T,\x)  =0 , \  1-\Delta \phi\geq 0\}$. In addition, any solution $\phi^*$ satisfies $1- \Delta \phi^* \geq (\frac{t}{T})^{\frac{d}{d+2}}$. 
\end{theorem}

While we will not repeat the proof here, let us mention that the proof starts with minimizing the Lyapunov (``entropy") functional among the very weak solutions $u$ of QPME 
\begin{equation}\label{Ljapunov}
    \int_{Q}  u^2(t,\x), 
\end{equation}
it can then be proved that the following formulations are  equivalent letting
\begin{itemize}
   
    \item $A: = \{u\in L^2(Q) \text{ is a very weak non-negative solution associated with }  u_0\in L^2 (\T)\} $, 
    \item $ B: = \{\phi \mid \phi(T,\x)  =0,\  1-\Delta \phi\geq 0
    \}$.
    
\end{itemize}
\begin{enumerate}
\item Original form
\begin{equation}
\begin{split}
    &I(u_0) = \inf_{u \in A}\sup_{\phi\in B} \int_{Q } \left(u^2 -2\partial_t\phi u -\Delta \phi u^2 +2u_0 \partial_t \phi
    \right)
\end{split}
\end{equation}
\item Flipping $\sup, \inf$
\begin{equation}\label{relaxed-form}
\begin{split}
    &J(u_0) = \sup_{\phi\in B} \inf_{u\in A} \int_{Q } \left(u^2 -2\partial_t\phi u -\Delta \phi u^2 +2u_0 \partial_t \phi
    \right)
\end{split}
\end{equation}
\item Point-wise minimization of \eqref{relaxed-form}.\\
\begin{equation}\label{phiformulation}
    \begin{split}
    &\tilde{J}(u_0)  =\sup_{\phi\in B}\  \int_{Q} \left(\frac{- (\partial_t \phi)^2}{1-\Delta \phi} + 2u_0\partial_t \phi
    \right) 
\end{split}
\end{equation}
\item Let $q = \partial_t \phi$, $\sigma  = 1-\Delta \phi$ in \eqref{phiformulation}
\begin{equation}\label{strong_terminal}
    \begin{split}
    &\hat{J} (u_0)  = \sup_{q, \sigma} \int_Q \left(
    \frac{-q^2}{\sigma} + 2 u_0 q
    \right)\\
     &\sigma \geq 0,\quad \sigma(T, \cdot) = 1,\quad  \partial_t \sigma+ \Delta \phi =0.
    \end{split}
\end{equation} 
\end{enumerate}
More specifically, it is proved that  
\begin{equation}\label{form_equivalency}
    \int_{Q}  u^2(t,\x) \ d\x\ dt = I(u_0) = J(u_0) = \tilde{J}(u_0) = \hat{J}(u_0).
\end{equation}

Theorem \ref{thm} shows that we can indirectly obtain very weak solutions to QPME by solving \eqref{phi_formulation}. We first obtain $\phi^{*}$, then obtain candidates for the very weak solution with \eqref{u_phi_thm}.
We can therefore consider the following loss function 
\begin{equation}\label{phi_form}
\begin{split}
\textbf{$\boldsymbol{\phi}$ formulation}\quad \mathcal{L}_{\phi}(\phi) &:= -\int_{Q} \left(\frac{- (\partial_t \phi)^2}{1-\Delta \phi} + 2u_0\partial_t \phi
    \right).
\end{split}
\end{equation}
It is not hard to see that if a smooth $\phi^*$ is a minimizer, as long as the recovered solution satisfies the homogeneous boundary condition and the initial condition, $\frac{\phi^*}{1-\Delta \phi^*}$ must be a solution to QPME. Moreover, in the case where $u\geq 0$, it has been proved that the solution to the QPME subject to $u_0\geq 0$ is unique. 
However, it is worth noting that such minimizer $\phi^*$ is not necessarily unique. Thus \eqref{phi_form} can be used to identify the unique solution to QPME, as long as the initial/boundary conditions are imposed, but more than one minimizer $\phi^*$ could theoretically exist \cite{brenier2020examples}.

Moreover, since \eqref{strong_terminal} is equivalent to \eqref{phiformulation}, one can also recover a candidate of very weak solution to QPME with $q^*$ and $\sigma^*$ by
\begin{equation}
    u_{q,\sigma}^*: = \frac{q^*}{\sigma^*}
\end{equation}
with $q^*$ and $\sigma^*$ being the maximizer of  \eqref{strong_terminal}. Thus, we may also consider the loss function 
\begin{equation}\label{qsigma_form}
  \textbf{$\mathbf{q}-\boldsymbol{\sigma}$ formulation}\quad  \mathcal{L}_{q,\sigma}(q,\sigma) = -\int_Q \left(\frac{-q^2}{\sigma} + 2 u_0 q\right).
\end{equation}

Similar to the discussion in section \ref{sec:PINN}, we can also relax the initial/boundary conditions of the recovered solution $u$ obtained from $\phi$ formulation and $q,\sigma$ formulation to a penalization by regularizing $\mathcal{L}_{B}$ and $\mathcal{L}_{I}$ as defined earlier:
\begin{equation}\label{full_phi}
    \begin{split}
        \mathcal{L}_{\phi-\text{NN}}(u): = \kappa  \mathcal{L}_{\phi}(u) + \mu\mathcal{L}_{B}(u) +
    \nu\mathcal{L}_{I}(u),
    \end{split}
\end{equation}
and 
\begin{equation}\label{partial_q_sigma}
    \mathcal{L}_{q,\sigma-\text{NN}}(u): = \kappa  \mathcal{L}_{q,\sigma}(u) + \mu\mathcal{L}_{B}(u) +
    \nu\mathcal{L}_{I}(u).
\end{equation}
However, it is worth pointing out that it is very difficult to impose initial condition as a hard constraint to solution ansatz for both formulations when the optimization problem is solved with neural networks as only intermediate functions $\phi,q,\sigma$ will be parametrized. The boundary conditions, on the other hand, can be explicitly imposed by modifying the solution ansatz.

For $q-\sigma$ formulation in particular, we should also note that, the consistency between $q$ and $\sigma$
\begin{equation}\label{qsigma_PDE}
    \partial_t \sigma + \Delta q = 0
\end{equation}
needs to be imposed since they are essentially derivatives of the same function $\phi$. This condition can be imposed by minimizing the residual of equation \eqref{qsigma_PDE} 
\begin{equation}
    \mathcal{L}_{\partial_t \sigma, \Delta q} =  \int_Q   ( \partial_t \sigma + \Delta q )^2 
\end{equation}
or in $L^1$ sense,
\begin{equation}\label{q_sigma_corelation}
    \mathcal{L}_{\partial_t \sigma, \Delta q} =  \int_Q   | \partial_t \sigma + \Delta q |.
\end{equation}
Thus, $\mathcal{L}_{q,\sigma -\text{NN}}$ can be further modified as 
\begin{equation}\label{qsigma_full}
    \mathcal{L}_{q,\sigma-\text{NN}} = \kappa\mathcal{L}_{q,\sigma}(u) + \mu\mathcal{L}_{B}(u) +
    \nu\mathcal{L}_{I}(u) +
    \gamma  \mathcal{L}_{\partial_t \sigma, \Delta q}.
\end{equation}

Let us remark that the condition \eqref{qsigma_PDE} can also be imposed weakly following the Dirichlet principle, $\mathcal{L}_{\partial_t \sigma, \Delta q}$ can be then replaced by

\begin{equation}
    \mathcal{L}_{\partial_t \sigma, \nabla q} :=\int_{[0,T]}\left(\frac{1}{2} \int_{\T}   |\nabla q|^2\ dx +  \frac{\lambda}{2} \left(\int_{\T} q\ dx\right)^2 + \int_{\T} \partial_t\sigma q\ dx \right)\ dt .
\end{equation}
While related results using $\mathcal{L}_{\partial_t \sigma, \nabla q}$ will not be presented in this paper, we remark that this formulation completely bypasses taking second order derivative of $q$, which means less smoothness requirement of the neural network ansatz. However, adding such a term would make the optimization more complicated especially in the high-dimensional cases. Since this term can not be interpreted as a pointwise condition as \eqref{q_sigma_corelation}, it thus can not benefit much from an efficient sampling scheme (see Section \eqref{sec:num_setting}).
In addition, the introduction of the extra hyper-parameter $\lambda$ further increases the difficulty of parameter tuning.  From our own experience, the training of the $q,\sigma$ formulation with the term $\mathcal{L}_{\partial_t \sigma, \nabla q}$ seems extremely challenging if not impossible and is therefore not presented.

\newcommand{\NN}{\mathcal{NN}}
\section{Solving high dimensional QPME with neural network ansatz}\label{sec:nn}
Neural network is a class of functions that have a certain layered structure, for example the feed-forward fully connected neural network is defined to be 
\begin{equation}\label{FFNN}
		\mathcal{NN}(\x; \theta) := W_n g(\cdots g(W_2 g(W_1\x+ b_1 )+b_2)\cdots) +b_n.
\end{equation}
In this case, each layer of the network is a composition of a linear transformation and an nonlinear function $g$ acting component-wise. Here, $\theta := [W_1, W_2,\cdots, W_n, b_1, b_2,\cdots, b_n]$ are the trainable parameters.

The idea of neural network based numerical solver for PDEs is to utilize such a neural network $\mathcal{NN}$ to approximate the function of interest, say $u$. This is usually achieved by solving an optimization problem 
\begin{equation}
    u   = \argmin_f \mathcal{C}(f),
\end{equation}
where $\mathcal{C}$ is some suitable objective function.  Then one could take a neural network as an ansatz and minimize $\mathcal{C}$ by tuning its parameters $\theta$ to get an approximate solution $\mathcal{NN}\left(\cdot; {\theta^*}\right)$ where
\begin{equation}
    \theta^* = \argmin_\theta \mathcal{C}\left(\mathcal{NN}\left(\cdot; {\theta}\right)\right).
\end{equation}
The process of optimization is also referred as ``training'', using the terminology from machine learning. The objective function $\mathcal{C}$ is often referred as the loss function.

In Section \ref{sec:formulation}, we have derived a few loss functions which can be used to solve the QPME. In this section, we provide further details on how to solve the aforementioned optimization problems with neural networks especially on how initial and boundary conditions can be imposed to the neural network as a solution ansatz. In particular, the following conditions are generally considered as a solution ansatz to QPME:
1) the initial condition, 2) the boundary condition, 3) the physical constraint, i.e., $u\geq 0$. In addition, one could consider to impose conditions like $1 - \Delta \phi \geq (\frac{t}{T})^{\frac{d}{d+2}}$ to narrow down the search space as we know by Theorem \ref{thm} that it is satisfied by the true solution.
We will slightly modify existing neural network structure as needed to satisfy the constraints.
In this paper, we take the architecture of the neural networks to be  feed-forward fully-connected as defined in \eqref{FFNN}, while other architectures could also be considered. 

\subsection{PINN formulation} 
To solve QPME with PINN formulation \eqref{PINN}, we first notice the argmin to \eqref{PINN_full} is exactly a solution to QPME. Then the solution itself can be directly parametrized with a neural network. In particular, to further impose the aforementioned conditions to the solution ansatz, we start with a neural network $\NNu(t,\x;\theta_u)$ with both time $t$ and spatial coordinates $\x$ as its inputs and denote the collection of trainable parameters as $\theta_u$. Moreover, since we need to compute the PDE residual, which 
includes the computation of second-order derivative of the solution ansatz, $\NNu(t,\x;\theta_u)$ must be at least second order differentiable. We thus require activation functions $g$ to be smooth ones such as $\tanh$ and $\text{softplus}$ functions.

To impose the \textbf{initial condition \eqref{IC} as a hard constraint}, we can parametrize the solution $u(t,\x)$ as:
$$u(t,\x;\theta_u) = u_0(\x) +t \nn_u(t,\x;\theta_u).$$
However, in this case, the physical constraint ($u \geq 0$) cannot easily be imposed explicitly. The positivity of the solution can only be reached through minimizing PDE residual. 

In the case where the \textbf{initial condition is imposed softly}, the term $\mathcal{L}_I$ defined as in \eqref{weak_initial} will be added as a part of the loss $\mathcal{L}_{\text{PINN}}$ and minimized through training. Meanwhile, the physical constraint of solution can be imposed by parametrization: 
$$ u(t,\x; \theta_u) = \text{softplus}\left( \nn_u \left(t,\x;\theta_u\right)\right)$$
where the softplus function is given by 
$$\text{softplus}(x) = \ln(1+e^x)$$
which guarantees the solution ansatz to be positive.

As for the boundary condition, 
the \textbf{homogeneous Dirichlet boundary condition \eqref{BC}} can be imposed as a hard constraint. We take advantage of the function  
\begin{equation}\label{f_dc}
    f_{dc}(\x) := \prod_{i=1}^{d} \frac{(a_i -x_i)(a_i+x_i)}{a_i^2}
\end{equation}
so that $f_{dc}(\x) = 0$ for any $\x\in\partial \T$.  Moreover, we notice that $f_{dc}(\mathbf{0}) = 1 $ and $0\leq f_{dc}(\x)\leq 1$ for all $\x\in \T$.
The solution ansatz $u(t,\x)$ can then be further modified by multiplying $f_{dc}$ to satisfy the boundary condition:  
\begin{equation}
    \begin{split}
    \textbf{hard I.C. + hard B.C.}\quad u(t,\x;\theta_u)  = u_0(\x)+tf_{dc}(x)\ \nn_u \left(t,\x;\theta_u\right)\\
        \end{split}
\end{equation}
\begin{equation}\label{PINN_with_condition}
    \textbf{soft I.C. + hard B.C.}\quad u(t,\x;\theta_u)  = f_{dc}(\x)\ \text{softplus} \left(\nn_u \left(t,\x;\theta_u\right)\right)
\end{equation}
assuming the homogeneous Dirichlet B.C. is satisfied by $u_0$.

The benefit of solving PINN is that the convergence of training of $\mathcal{L_{\text{PINN}}}$ guarantees accurate solution, which is justified in \eqref{convergence_L1} and \eqref{convergence_L2}.
On the other hand, PINN formulation also has its own limitation that it only allows strong solutions. Solutions with less regularity can not be identified with this formulation.

\subsection{\texorpdfstring{$\phi$}{} formulation}

To solve QPME following the $\phi$ formulation, we need to parametrize $\phi(t,\x)$ in \eqref{phi_form} instead of the solution $u$ directly.
When computing $\mathcal{L}_{\phi}$ as in \eqref{phi_form}, we also need the ansatz of $\phi(t,\x)$ to be at least second-order differentiable. Note that this is a much weaker assumption on the solution ansatz of $u$ compared with the PINN. In particular, no assumption is needed on the smoothness of $u$ directly. We simply take a neural network $\nn_{\phi} (t,\x;\theta_{\phi})$ with smooth activation function as its solution ansatz.

We then note that the minimizers $\phi^*$ to \eqref{phi_form} must also satisfy certain conditions in order to obtain reasonable solutions.
As suggested by Theorem \ref{thm}, we would like to require $\phi$ to vanish at $t = T$. We thus let 
 \begin{equation}
     \phi(t,\x;\theta_{\phi}) = (T-t)\nn_{\phi}(t,\x;\theta_{\phi}).
 \end{equation}

For the recovered solution $u_{\phi}$
 \begin{equation}\label{u_phi}
     u_{\phi}:= \frac{\partial_t \phi }{1-\Delta \phi},
 \end{equation} 
 unlike PINN formulations, the solution to QPME is not directly parametrized; thus it is not easy to impose the initial condition as a hard constraint. Instead, we enforce the constraint softly relying on the penalty term $\mathcal{L}_{I}$.
 The homogeneous Dirichlet boundary condition, on the other hand, can be softly enforced with the term $\mathcal{L}_B$ or  enforced as a hard constraint by modifying the neural network. Essentially, we only need 
 $$\partial_t \phi|_{\partial \T} = 0,$$
 which can be achieved using the ansatz 
 \begin{equation}\label{phi_condition}
    \textbf{soft I.C. + hard B.C.} \quad \phi (t,\x;\theta_{\phi})  = (T-t) f_{dc}(\x) \nn_{\phi} (t,\x;\theta_{\phi}),
 \end{equation}
 where $f_{dc}(\x)$ is defined as in \eqref{f_dc}. 
Additionally, while the recovered solution $u_{\phi}\geq 0$ is desired, this condition cannot be easily imposed  by simply modifying the solution ansatz.

Compared with the PINN formulation, using a $\phi$-formulation allows solutions in a very weak sense. It potentially can find solutions with less regularity. The smoothness requirement is not directly applied to $u_{\phi}$. However, as described above, a few conditions of $\phi$ can not be easily enforced. In addition to the positivity of $u_{\phi}$, conditions like 
 \begin{equation}\label{growing_condition}
      1- \Delta \phi \geq (\frac{t}{T})^{\frac{d}{d+2}},
\end{equation}
is difficult to enforce as a hard constraint either. While \eqref{growing_condition} is preferable as it can narrow down the search function space for $\phi$ (since we know the PDE solution would satisfy that), it is not necessary. However, the fact that $u_{\phi}$ is not confined to be positive function can potentially cause the training of $\mathcal{L}_{\phi-NN}$ converges to unphysical solutions. 

\subsection{\texorpdfstring{$q$}{}-\texorpdfstring{$\sigma$}{} formulation} 
The $q-\sigma$ formulation \eqref{qsigma_form} is derived from the $\phi$ formulation, and thus also inherits a few conditions for $\phi$.
We first notice that when computing $\mathcal{L}_{q,\sigma}$, no computations of derivatives will be needed. However, when computing $\mathcal{L}_{\Delta q,\partial_t\sigma}$, first-order and second-order derivatives of $\sigma$ and $q$ are required respectively. We thus can start by parametrizing $q(t,\x)$ and $\sigma(t,\x)$ with neural networks $\nn_q(t,\x;\theta_{q})$ and $\nn_{\sigma}(t,\x;\theta_{\sigma})$ which should be at least first and second order differentiable respectively. 

To ensure the positivity of $\sigma$,
as suggested in \eqref{strong_terminal}, we further parametrize $\sigma$ by
$$\sigma (t,\x ;\theta_{\sigma})= \text{softplus}\left( \nn_{\sigma}\left(t,\x ;\theta_{\sigma}\right)\right).$$
To guarantee that 
$$\sigma(T,\cdot) = 1,$$
we modify the above and let
\begin{equation}\label{sigma_condition}
    \sigma(t,\x;\theta_{\sigma})  = \text{softplus}\big(\ln\left(e-1\right)+\left(T- t\right) \nn_{\sigma}\left(t,\x;\theta_{\sigma}\right) \big).
\end{equation}
Alternatively, if we also impose the condition (to narrow down the search space)
\begin{equation}\label{sigma_selection}
     \sigma \geq (\frac{t}{T})^{\frac{d}{d+2}},
\end{equation}
one can also parametrize $\sigma$ as 
\begin{equation}
    \sigma (t,\x;\theta_{\sigma}) = \left(\frac{t}{T}\right)^{\frac{d}{d+2}}+(T- t)\, \text{softplus}\left(\nn\left(t,\x;\theta_{\sigma}\right) \right).
\end{equation}

For the recovered solution 
$$u_{q,\sigma} = \frac{q}{\sigma},$$ like $\phi$ formulation, the initial condition can only be softly imposed. The homogeneous Dirichlet boundary condition can be enforced as a hard constraint, as long as 
$$q|_{\partial \T} = 0.$$ We thus let 
\begin{equation}
   q(t,\x;\theta_{q}) = f_{dc}(\x)\nn_{q}(t,\x;\theta_{q}).
\end{equation}
Moreover, to ensure $u_{q,\sigma} \geq 0$, we further let 
\begin{equation}\label{q_with_conditions}
   \textbf{soft I.C. + hard B.C. }\quad q(t,\x;\theta_{q}) = f_{dc}(\x)\text{softplus}\left(\nn_{q}\left(t,\x;\theta_{q}\right)\right).
\end{equation}
Similar to $\phi$ formulation, the $q-\sigma$ formulation allows solutions with less regularity. However, two neural networks  will be needed to parametrize a solution to QPME, which could potentially be more challenging to train.

\subsection{Empirical loss and training data sampling}\label{empirical_loss}
\newcommand{\p}{\mathrm{P}}
\newcommand{\PINN}{\mathcal{L}_{\text{PINN}}}
\newcommand{\X}{\mathbf{X}}
To solve the QPME with aforementioned formulations \eqref{PINN_full},\eqref{full_phi} and \eqref{qsigma_full}, we need to compute the high-dimensional integrals of the neural network or its derivatives to evaluate the loss functions. 

In practice, Monte Carlo methods are usually used to approximate those high dimensional integrals. The approximate solutions are then obtained by minimizing the surrogate empirical loss functions. Take the PINN formulation as an example, let $\p_{\T}$ be the uniform probability distributions over the spatial domain $\T$ and let $\{\X_j\}_{j=1}^{n}$ be an i.i.d. sequence of random variables distributed according to $\p_{\T}$. Parallelly, we also define $\p_{[0,T]}$ be the uniform probability distributions over the spatial domain $[0,T]$ and let $\{T_j\}_{j=1}^{n}$ be an i.i.d. sequence of random variables distributed according to $\p_{[0,T]}$.  Define the empirical loss $\PINN^{n}$ by setting
\begin{equation}\label{empirical_PINN}
  \PINN^{n} =   \frac{ \kappa}{n}\sum_{j=1}^{n}\left( \partial _t u(T_j,\X_j) -\frac{1}{2} \Delta u(T_j,\X_j)\right)^2 + \frac{\nu}{n} \sum_{j=1}^{n}\left(u \left(0,\X_j\right) - u_{0}\left(\X_j\right)\right)^2
\end{equation}
for the case where only I.C. is imposed softly and the loss measuring norm is taken to be $L^2$. Notice that all terms are scaled by $\frac{1}{|\T|}$, which does not change the minimizer of the problem but can effectively avoid numerical blowup in evaluating the loss during the training. Similarly, $\mathcal{L}_{B}$ can also be approximated with points uniformly sampled from $\partial \T$ when needed. We further refer such sampled data as training data.

However, a uniform sampling of $\X_j$'s sometimes does not meet the need of our computation especially in the case where the dimension $d$ is very large. Notice that one essential feature of solutions to QPME is that it has a free boundary that separates positive part of the solution from the zeros. In particular, in the case where the solution is a Barenblatt solution, the nonzero values of the solution actually concentrate near the origin. Ideally, one would like to sample points in both non-zero and zero region, to capture the local features of the solution. However, with a fixed budget of training samples, it could happen that all randomly sampled data points only reside in the zero region, which is apparently problematic. In fact, this could become a serious issue for high dimensional problems. For example, when $d =20$, the probability of sampling the nonzero region of a Barrenblatt solution \eqref{barenblatt} at $t=2$ within $[-7,7]^{20}$ can be computed by the ratio of the volume of the $d$-ball $V_{\text{nonzero}}$ with radius $(22)^{1/2} 2^{6/11}$  standing for the non-zero region versus the volume of the hypercubic. It can then be computed that 
$$ \p_{\text{nonzero}} = \frac{V_{\text{nonzero}}}{14^{20}}\approx 1.57\times 10^{-8}$$
Which means, the non-zero region can rarely be sampled if not impossible. 
Therefore, we would need an effective sampling scheme which puts larger weights over the non-zero region, so that we can accurately approximate the loss function.  Ideally, one could hope for an adaptive important sampling scheme which provides samples of the training data based on a distribution adapted to current status of the parameterized solution and its derivatives throughout the training process. However, especially for high-dimensional problems, such sampling scheme is challenging and computationally expensive to implement. Therefore, a hand-crafted sampling scheme is used instead, which is explained in details in Section \ref{sec:nuemrical} for specified numerical examples.

\newcommand{\tX}{\tilde{\X}}

\section{Numerical example: Barenblatt Solution }\label{sec:nuemrical}
To test the numerical schemes, we will use a series of special solutions, known as
Barenblatt Solutions. They are given by \begin{equation}
    U_m(t, \x; C) := t^{-\alpha}\left(\left(C - \frac{\beta(m- 1)}{2}\frac{|\x|^2}{t^{2\beta} }\right)^+\right)^{\frac{1}{m-1}}
\end{equation}
where $\alpha := \frac{d}{d(m-1)+2}, \beta := \frac{\alpha}{d}$, $(s)^+ := \max(s,0)$ and $C > 0$ is an arbitrary constant. This solution takes a Dirac mass as initial data : $\lim_{t \to 0} u(t,\x) \to M\delta(\x)$, where $M$ is a function of the constant $C$ (depends on $m$ and $d$).
In the particular case $m = 2$, the Barenblatt Solution to \eqref{QPME} reduces to  
\begin{equation}\label{barenblatt}
    U_2(t, \x; C) := t^{-\frac{d}{d+2}}\left(C - \frac{1}{2(d+2)}\frac{|\x|^2}{t^{\frac{2}{d+2}} }\right)^+.
\end{equation}
The free boundary $\partial \Pu$ to \eqref{barenblatt} in this case can then characterized by the equation 
$$ \quad|\x| = r_t$$
with $r_t:=  \sqrt{2C(2+d)}\ t^{\frac{1}{d+2}}$. We also notice that the solution is scale-invariant:
    \begin{equation*}
    u_{\lambda}(t, \x) := \lambda^{\alpha} u(\lambda t, \lambda^{\beta} \x).
    \end{equation*}
The shifted Barenblatt Solution is a {strong/weak/very weak} solution of PME, which is unique subject to the Dirac initial condition. 

Since numerically one cannot set $\delta$ function as the initial condition, we specifically consider the following IBVP:
\begin{equation}\label{numerical_baren}
\begin{split}
     &\partial_t u =  \frac{1}{2}\Delta u^2\quad (t,\x)\in Q = [0,1]\times\T,\\
     &u(0,\x) = \left(1- \frac{1}{2(2+d)} |\x|^2
     \right)^+.
\end{split}
\end{equation}
Notice the initial condition is essentially the Barenblatt Solution \eqref{barenblatt} evaluated at $t=1$ when $C =1$.
The exact solution to \eqref{numerical_baren} is therefore the Barenblatt Solution \eqref{barenblatt} with the time shifted: 
 \begin{equation}\label{exact}
    U_2(t, \x) := \left(t+1\right)^{-\frac{d}{2+d}}\left(1 - \frac{1}{2(d+2)}\frac{|\x|^2}{ (t +1)^{\frac{2}{d+2}} }\right)^+.
\end{equation}
We further let $\T = [-a,a]^d$, where $a$ is the smallest integer that is greater than the radius of the free boundary of $U_2(t,\x)$ at the terminal time $T=1$: 
$$a := \text{ceil}(r_T)$$
where $r_T  = (2+d)^{\frac{1}{2}}2^{\frac{4+d}{2d+4}}$ to ensure the the computational domain is large enough to include the entire free boundary for $t\in [0,1]$.
We take this example to test the effectiveness of the proposed formulations by comparing the approximate solution with the exact one \eqref{exact}. The performance of each formulation is further analyzed to show the pros and cons.

\subsection{Numerical settings}\label{sec:num_setting}
In particular, We would like to solve the aforementioned QPME with three formulations using neural network ansatz.
We specifically take $\nn_u(\cdot,\cdot;\theta_u)$, $\nn_{\phi}(\cdot,\cdot;\theta_{\phi})$, $\nn_{q}(\cdot,\cdot;\theta_{q})$ and $\nn_{\sigma}(\cdot,\cdot;\theta_{\sigma})$ to be fully connected neural networks with two hidden layers and the $\text{softplus}(\cdot)$ as their activation function. The corresponding solution ansatz can then be developed following Section \ref{sec:nn}. Notice that when computing the derivatives of the solution ansatz, a chain rule will be applied and the derivatives of $f_{dc}(\x)$ needs to be computed when a homogeneous Dirichlet B.C. is imposed.

To evaluate the empirical losses, as discussed in Section \eqref{empirical_loss}, we further take randomly sampled data to approximate the integrals in $Q$ or in $\T$. When boundary condition is softly imposed, we take additional randomly sampled data over $\Sigma_T$ to evaluate $\mathcal{L}_B$. In particular, since the data are sampled on the fly, a new set of data will be sampled at each training step, the total number of training data $n$ can then be computed by $n= \text{batch size}\times  \text{training steps}$.

\begin{figure}[h!]
    \centering
\begin{subfigure}[t]{0.45\textwidth}
  \includegraphics[width = 1.2\textwidth]{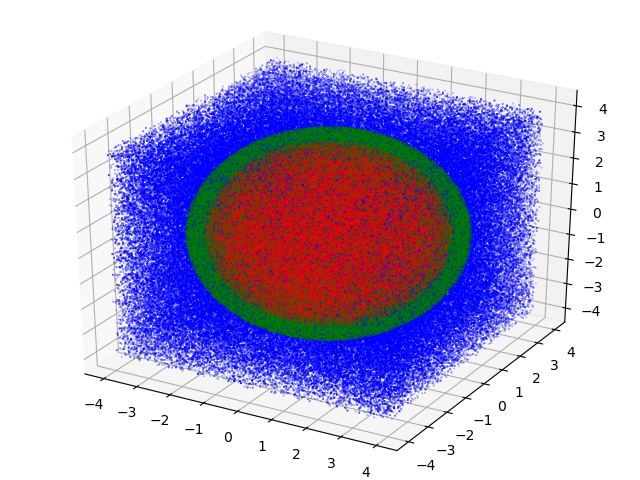}
        \caption{\textbf{3d:} $\theta_0 = 0.3, \theta_0 = 0.3$, $\T = [-4,4]^3$.}

\end{subfigure}
    \hfill
\begin{subfigure}[t]{0.45\textwidth}
      \includegraphics[width = 1.2\textwidth]{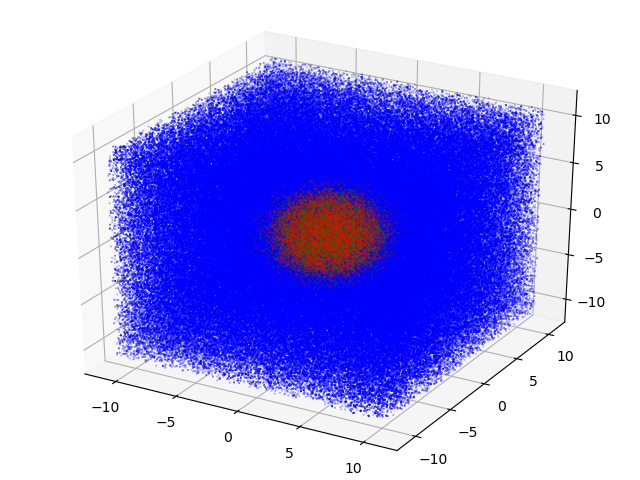}
      \caption{\textbf{50d:} $\theta_0 = 0.3, \theta_1 = 0.2$, $\T = [-11,11]^{50}$ }
    \end{subfigure}
    \caption{3D projection (first three coordinates) of samples of $\{\tX_j\}_{j=1}^{10^6}$ in $\T$. Red : $\tX_j \in V_0$, green: $\tX_j \in V_1$ and blue: $\tX_j \in V_2$.}
    \label{fig:train_data}
\end{figure}

Moreover, for high-dimensional cases, to make sure the sampled points can capture the features of the solutions, we use the following weighted sampling scheme. Specifically, we first decompose the region $\T = [-a,a]^d$ into $\T =  V_0\cup V_1\cup V_2$ where $$V_{0}: =\{\x \in \T\ |\ |\x|\leq r_0 \},\quad V_{1}: =\{\x \in \T\ |\ r_0<|\x|\leq r_T \},\quad V_{2}: =\{\x \in \T\ |\ |\x|> r_T \}.$$
The radius of these region are decided by the radius of the free boundary to the Barenblat solution \eqref{exact} at $t=0$ and $t=T=1 $ respectively: 
$$r_0 : = \sqrt{2(2+d)},\qquad r_T  = (2+d)^{\frac{1}{2}}2^{\frac{4+d}{2d+4}}.$$
We then take  weights $\theta_0, \theta_1$ and 
the $\tilde{\X}_j$'s will be uniformly sampled within $V_0$ with a probability $\theta_0$, within $V_1$ with a probability $\theta_1$ and within $\T$ with a probability $\theta_2:=1-\theta_0-\theta_1$ (see Figure \ref{fig:train_data} for an illustration of the sampled training data).  
The probability of sampling data in each region can be thus be computed as 
$$P_{V_0} = \theta_0+\theta_2 \frac{|V_0|}{|\T|}\quad P_{V_1} = \theta_1+\theta_2 \frac{|V_1|}{|\T|} \quad P_{V_2}=\theta_2 \frac{|V_2|}{|\T|}.$$
The density function to this mixture distribution can be written as the piecewise constant function
$$f(\x) = \theta_0 f_0(\x)+ \theta_1 f_1(\x)+\theta_2 f_2(\x)$$
where 
\begin{equation*}
    \begin{split}
        f_0(\x) = \frac{1}{|V_0|}\mathbf{1}_{V_0}(\x), \quad
        f_1(\x) = \frac{1}{|V_1|}\mathbf{1}_{V_1}(\x), \quad
        &f_2(\x) =\frac{1}{|\T|}
    \end{split}
\end{equation*}
are defined over the entire $\T$,  with $\mathbf{1}_{V}(\x)$ being an indicator function of region $V$.

When sampling from $V_0$ and $V_1$, to ensure the data are uniformly sampled in those high-dimensional ball regions, we specifically adopt the following algorithm as in \cite{marsaglia1972choosing}
\begin{enumerate}
    \item Generating random points uniformly on the $(d-1)$-unit sphere
    \begin{enumerate}
        \item Generate an $d$-dimensional vector $\bm{x} = (x_1x_2,\cdots, x_d)$ so that $x_i \sim N(0,1)$ for $i= 1,2,\cdots, d$.
        \item then $\tilde{\bm{x}} := \frac{\bm{x}}{||\bm{x}||_2}$ is a uniformly sampled point form $(d-1)$-unit sphere.
    \end{enumerate}
    \item Generate a point uniformly at random {\it{within}} the $d$-ball
    \begin{enumerate}
        \item  Let $u$ be a number generated uniformly at random from the interval $[0, 1]$, then $u^{\frac{1}{d}}\tilde{\bm{x}}$ is a point randomly sampled within the unit ball.
        \item Further, $r_0 u^{1/d}\tilde{\bm{x}} $ is a random point in $V_0$ and  $\left((r_T-r_0) u^{1/d}+r_0\right)\tilde{\bm{x}} $ is a random point in $V_1$.
    \end{enumerate}
\end{enumerate}

To avoid changing the values of the integrals to be evaluated in the loss functions, a piecewise constant factor should be multiplied to the empirical function to correct the approximation of the integrals resulted from a nonuniform distributed training data. For PINN formulation, the empirical loss \eqref{empirical_PINN} can then be rewritten as 
\begin{equation}
  \PINN^{n} =   \frac{\kappa}{n}\sum_{j=1}^{n} c(\tX_j)\left( \partial _t u(T_j,\tX_j) -\frac{1}{2} \Delta u(T_j,\tX_j)\right)^2 +\frac{\nu}{n}\sum_{j=1}^{n} c(\tX_j)\left(u \left(0,\tX_j\right) - u_{0}\left(\tX_j\right)\right)^2
\end{equation}
with the correction term
$$c(\x) = \displaystyle \sum_{i=1}^{3}\frac{|V_i|}{|\T| P_{V_i}}\mathbf{1}_{V_i}(\x).$$
Here, $\tX_j$'s are random points in $\T$ sampled subject to the aforementioned density function $f(\x)$, while $T_j$ will be uniformly sampled from $[0,1]$. The empirical losses for $\phi$ formulation and  $q-\sigma$ formulation can also be formulated in a similar fashion.
The data are then randomly sampled on the fly which are batched into $1000$ for each training step. New data is sampled for each batch. Essentially, this sampling scheme enforces data samples in all three regions which can potentially improve the representativity of the training data and thus lead to a faster convergence of the training procedure.

For high-dimensional cases, the initial conditions and PDEs are in fact imposed without the correction factor $c(\x)$, i.e., with the efficient sampling scheme, we are essentially minimizing the modified initial/PDE conditions. Taking such terms measured with $L^2$ as examples, the following loss terms will be minimized
\begin{equation}
\begin{split}
    &\mathcal{L}_{I}(u) = \int_{\Omega} \left(u(0,\x)-u_0\left(\x\right)\right)^2 f(\x)\ d\x, \\ 
    &\mathcal{L}_{\text{PDE}} (u) = \int_Q \left( \partial _t u -\frac{1}{2} \Delta u^2\right)^2 f(\x)\ d\x dt,\\
    &\mathcal{L}_{\partial_t \sigma, \Delta q} =  \int_Q   ( \partial_t \sigma + \Delta q )^2 f(\x)\ d\x dt.
\end{split}
\end{equation}
Since $f(\x)$ is merely a positive piecewise constant function, such modification will keep the minimizers to these terms unchanged, meaning the desired initial condition and PDE equation will still be imposed with mild assumption on the regularity of the solution ansatz.
The reason that the correction constants are not used is because for high-dimensional cases, $$c(\x)\ll 1\% \quad\forall \x \in V_0\cup V_1 ,$$
which means extremely small contribution of samples within these region will be made to update the trainable parameters with SGD. 
We also notice that while this is possible for imposing the initial condition and PDEs, we must apply $c(\x)$ the to the inf terms, namely $\mathcal{L}_{\phi}(\phi)$ and $\mathcal{L}_{q,\sigma}(q,\sigma)$ as the sampling scheme will  change the optimization target.

However, the choices of $\theta_0$ and $\theta_1$ still remain to be arbitrary. While numerical examples show that certain choices could lead to faster convergence, there is no clear principles one could follow to make optimal choices. 
Similar situations happen when we decide the values of $\nu, \kappa,\gamma$ to balance the terms in the losses. While theoretically these hyper-parameters can be any positive number, the choice of them can heavily influence the training procedure. Some choices seem to help the weighted loss to converge faster compared to others, but there is no justified reason for any certain choices. Therefore, the choices of theses hyper-parameters used in the results reported in this paper are results of trial and error.

The losses are then optimized by tuning the trainable parameters of the neural networks. We take the optimizer that implements the Adam algorithm \cite{kingma2014adam} to train the models. The complete algorithm to establish the loss function and to train the neural networks is implemented using the TensorFlow library~\cite{abadi2016tensorflow}.

Once the training is finished, to evaluate the quality of the approximate solution obtained with the trained neural networks, we can further quantify the generalization error of it. In particular, we define the relative errors on a solution slice $u(t,x,y,c,\cdots,c)$ at time $t$ for some fixed constant $c\in[-a,a]$, denoting $u(t,x,y,c,\cdots,c)$ by $u(t)$ for simplicity:
\begin{equation}
\begin{split}
    & L^1\textbf{-Relative Error}\quad \frac{||u_{NN}(t)-u(t)||_{1}}{||u(t)||_{1}},\\
    & L^2\textbf{-Relative Error}\quad \frac{||u_{NN}(t)-u(t)||_{2}}{||u(t)||_2},\\
   & H^1\textbf{-Relative Error}\quad \frac{||u_{NN}(t)-u(t)||_{H^1}}{||u(t)||_{H^1}}.\\
\end{split}
\end{equation}
where $u_{NN}$ stands for the neural network based solutions.
These norms can be further numerically approximated over a $100\times 100$ evenspaced mesh on $[-a,a]^2$ letting $\{x_i,y_j\}_{i=1,j=1}^{100}$ be the meshgrid points
\begin{equation*}
    \begin{split}
        &||f(x,y)||_1 \approx \frac{(2a)^2}{10^4}\sum_{i,j}|f(x_i,y_j)|,\\
        &||f(x,y)||_2 \approx \sqrt{\frac{(2a)^2}{10^4}\sum_{i,j}|f(x_i,y_j)|^2},\\
        &||f(x,y)||_{H^1} \approx \sqrt{\frac{(2a)^2}{10^4}\sum_{i,j}(|f(x_i,y_j)|^2 + |\nabla f(x_i,y_j)|^2)}.\\
 \end{split}
\end{equation*}
The numerical relative errors can then be computed with predicted values of the neural network solutions evaluated at the mesh-grid points.

\subsection{PINN formulation}
In this section, we  consider the case where the QPME \eqref{numerical_baren} is solved following a PINN formulation \eqref{PINN_full} in both $L^1$ and $L^2$ norm. Specifically speaking, the homogeneous Dirichlet boundary condition is imposed as a hard constraint and the initial condition is imposed as a soft constraint following \eqref{PINN_with_condition}. The specific algorithmic settings are further presented in Table \ref{tab:PINN_error} along with the relative errors computed for the trained solution slice $u(0.5,x,y,1.0,\cdots,1.0; \theta_u^*)$ at time $t= 0.5$ comparing with the exact solution \eqref{exact}.

\begin{figure}[htbp]
     \centering
     \begin{subfigure}[t]{0.3\textwidth}
         \centering
         \includegraphics[width=1.2\textwidth]{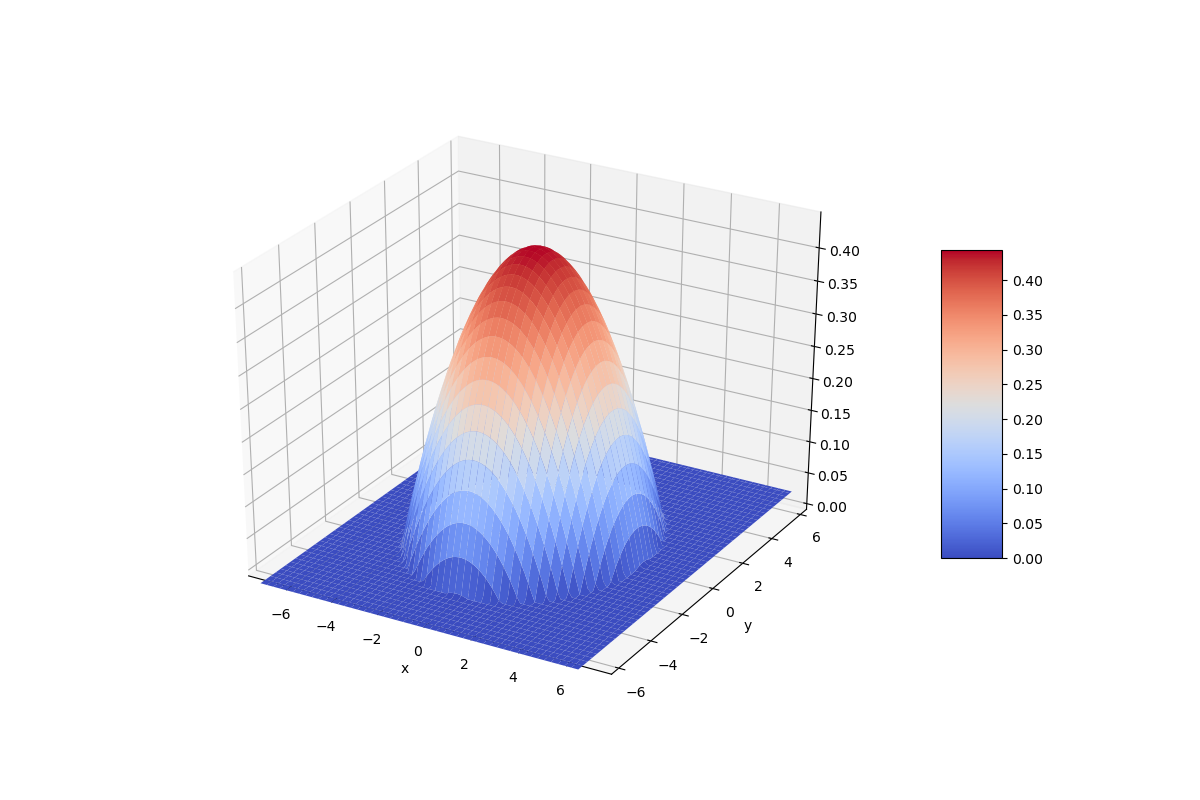}
         \caption{Barenblatt reference solution  }
     \end{subfigure}
    \hfill
     \begin{subfigure}[t]{0.3\textwidth}
         \centering
         \includegraphics[width=1.2\textwidth]{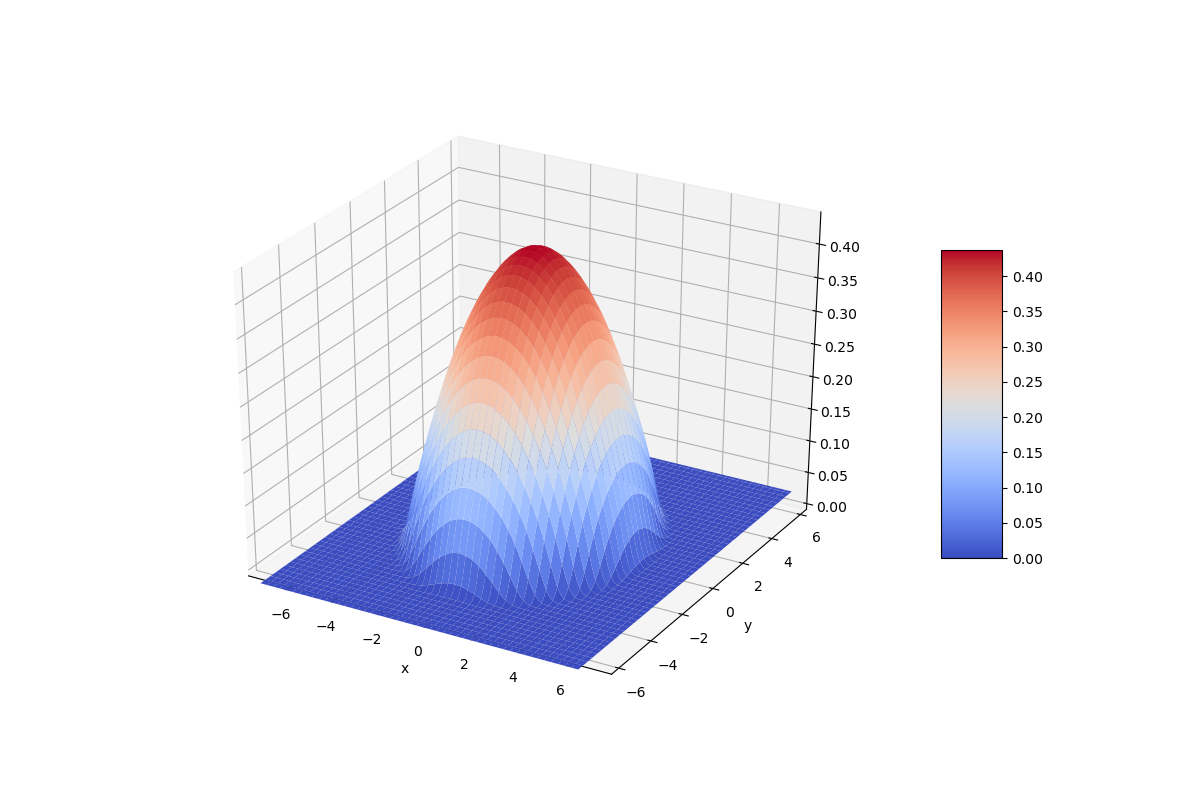}
         \caption{Learned solution slice}
         
     \end{subfigure}
     \hfill
    \begin{subfigure}[t]{0.3\textwidth}
         \centering
         \includegraphics[width=1.2\textwidth]{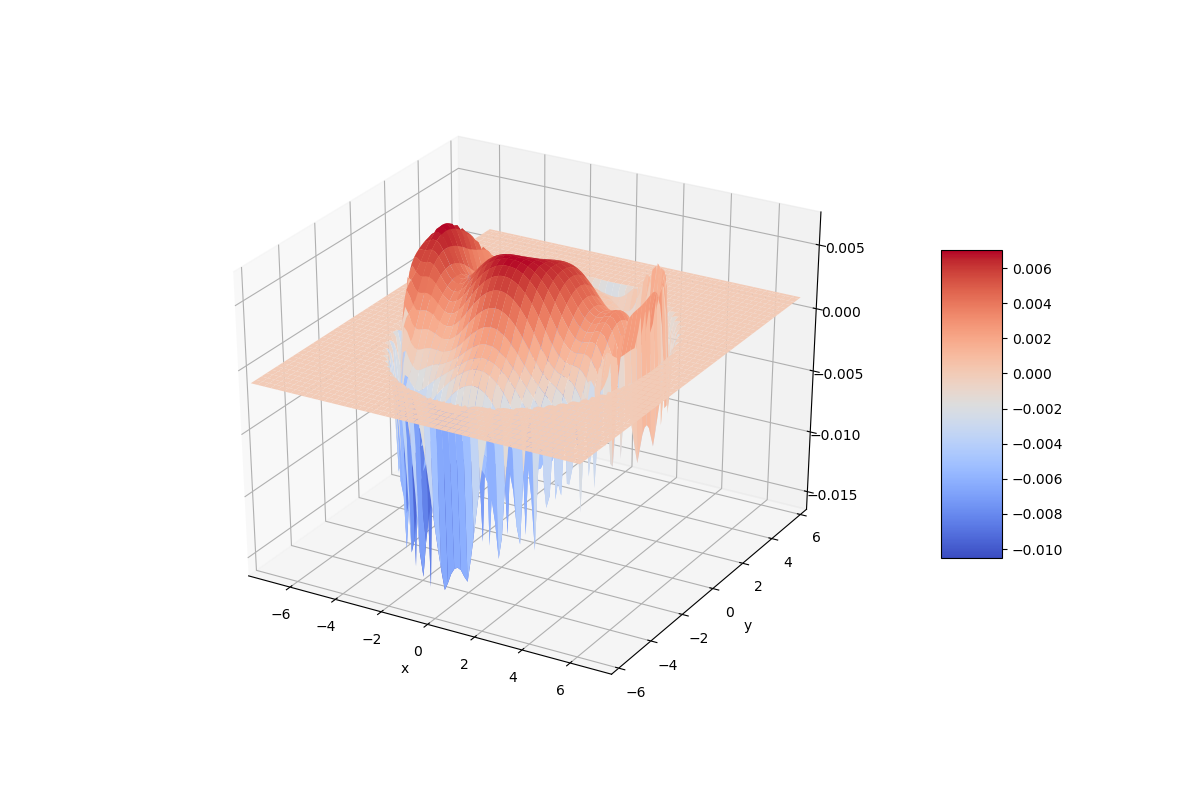}
         \caption{Learned solution error}
         
     \end{subfigure}
    \hfill
    \begin{subfigure}[t]{0.3\textwidth}
         \centering
         \includegraphics[width=1.2\textwidth]{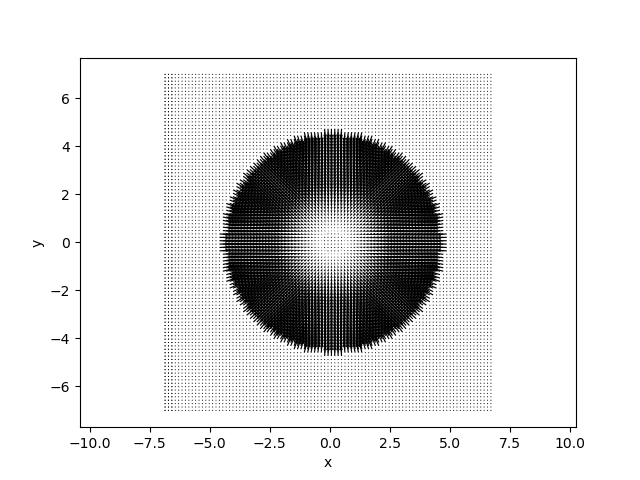}
         \caption{Barenblatt reference solution gradient}
         
     \end{subfigure}
     \hfill
    \begin{subfigure}[t]{0.3\textwidth}
         \centering
         \includegraphics[width=1.2\textwidth]{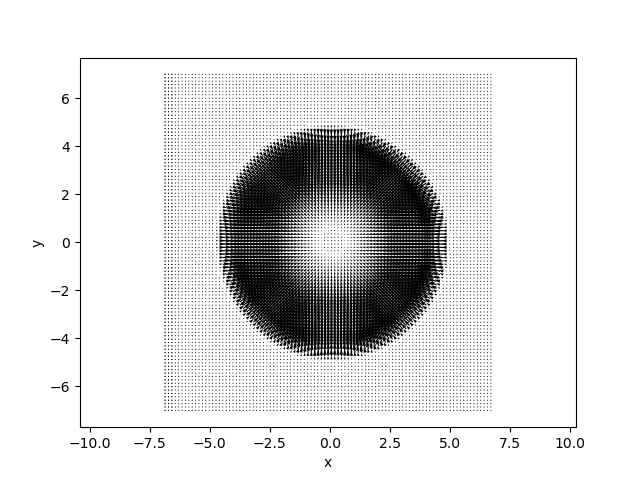}
         \caption{Learned solution gradient}
         
     \end{subfigure}
      \hfill
       \begin{subfigure}[t]{0.3\textwidth}
         \centering
         \includegraphics[width=1.2\textwidth]{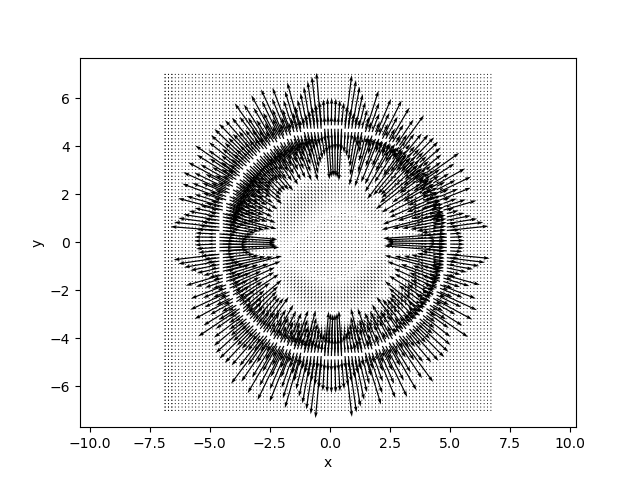}
         \caption{Learned solution gradient error}
         
     \end{subfigure}
        \caption{\textbf{15D, $L^2-$ PINN formulation \eqref{PINN_full}} Predicted solution slice $u(0.5,x,y,1.0,\cdots, 1.0)$ for $\x\in \T = [-7,7]^{15}$, $t= 0.5$. }
         \label{fig:PINN_l2}
\end{figure}

\begin{figure}
    \centering
    
    \begin{subfigure}[t]{0.47\textwidth}
             \centering
    \includegraphics[width = 1.2\textwidth]{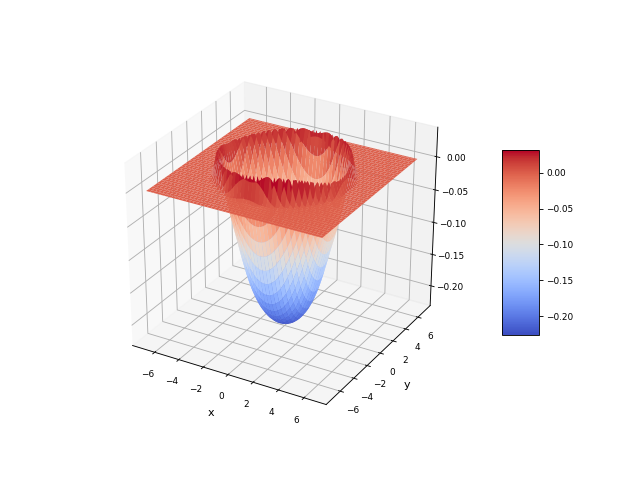}
               \caption{Learned $u_t$}
    \end{subfigure}
    \hfill
    \begin{subfigure}[t]{0.47\textwidth}
             \centering
    \includegraphics[width = 1.2\textwidth]{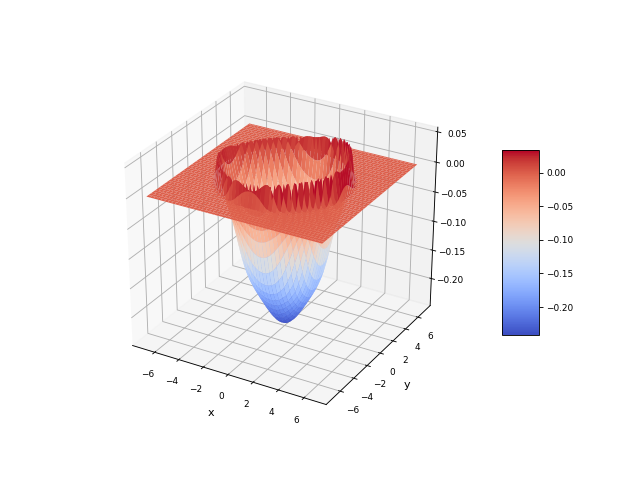}
               \caption{Learned $\displaystyle \frac{1}{2}\Delta u^2$}
    \end{subfigure}
    
    \caption{\textbf{15D, $L^2-$ PINN formulation \eqref{PINN_full}},  predicted partial derivatives.}
    \label{fig:PDE_l2}
\end{figure}

\begin{figure}
    \centering
    \begin{subfigure}[t]{0.47\textwidth}
             \centering
    \includegraphics[width = 1.2\textwidth]{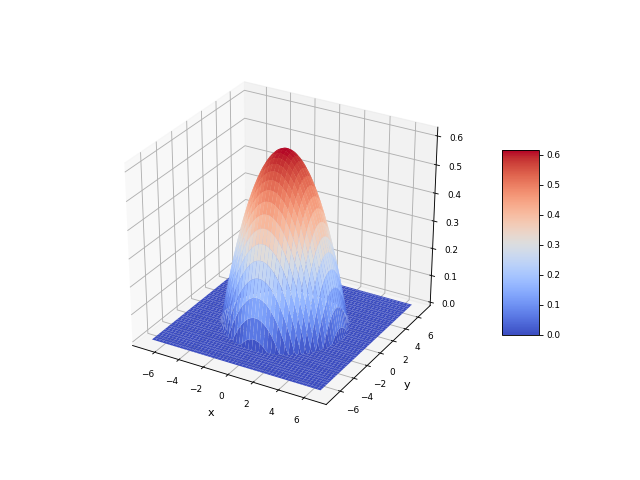}
               \caption{Exact $u_0$}
    \end{subfigure}
    \hfill
    \begin{subfigure}[t]{0.47\textwidth}
             \centering
    \includegraphics[width = 1.2\textwidth]{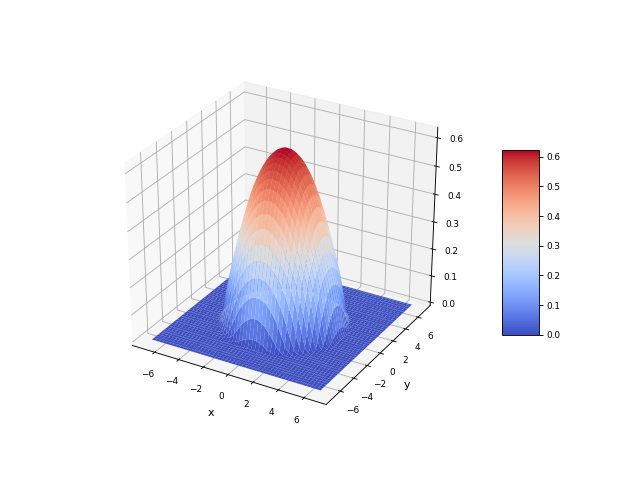}
               \caption{Learned initial value}
    \end{subfigure}
    
    \caption{\textbf{15D, $L^2-$ PINN formulation \eqref{PINN_full}}, predicted initial value $u(0,x,y,1.0,\cdots,1.0)$ for $\x \in 
    \T =[-7,7]^{15}$.}
    \label{fig:init_l2}
\end{figure}

\begin{figure}[htbp]
     \centering
     \begin{subfigure}[t]{0.3\textwidth}
         \centering
         \includegraphics[width=1.2\textwidth]{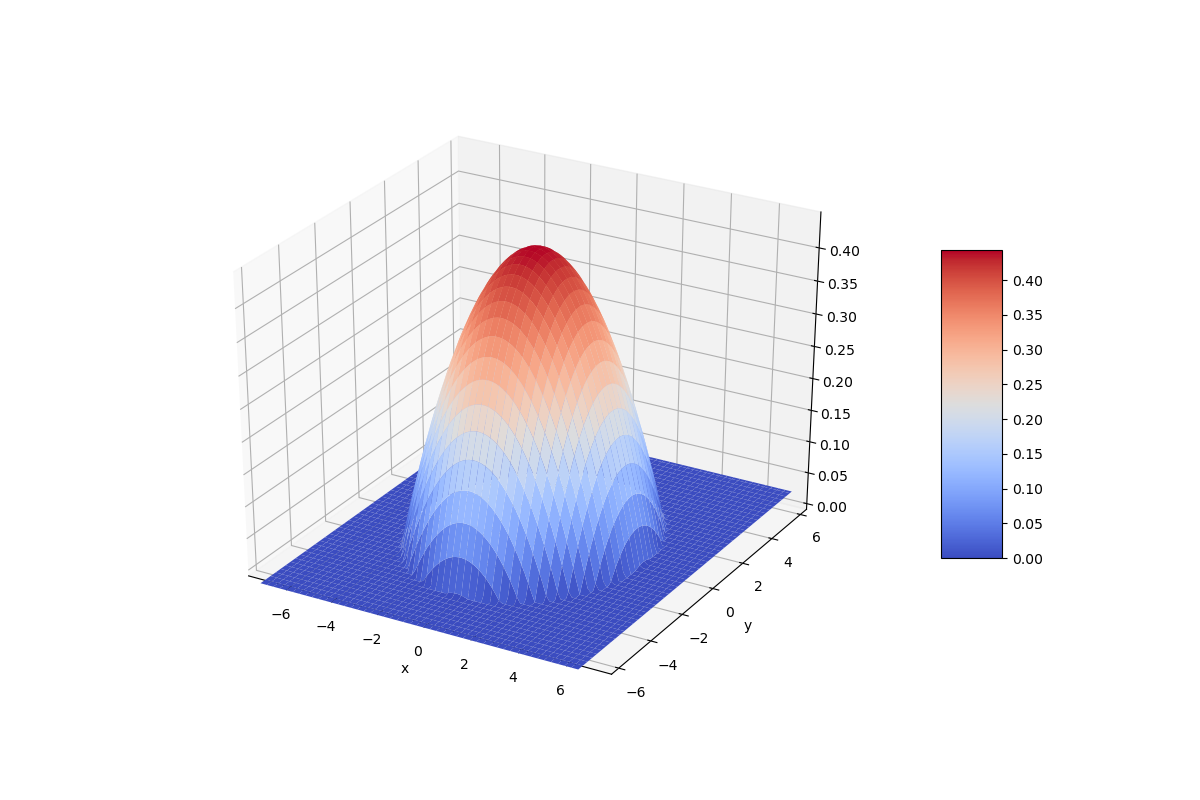}
         \caption{Barenblatt reference solution  }
         
     \end{subfigure}
    \hfill
     \begin{subfigure}[t]{0.3\textwidth}
         \centering
         \includegraphics[width=1.2\textwidth]{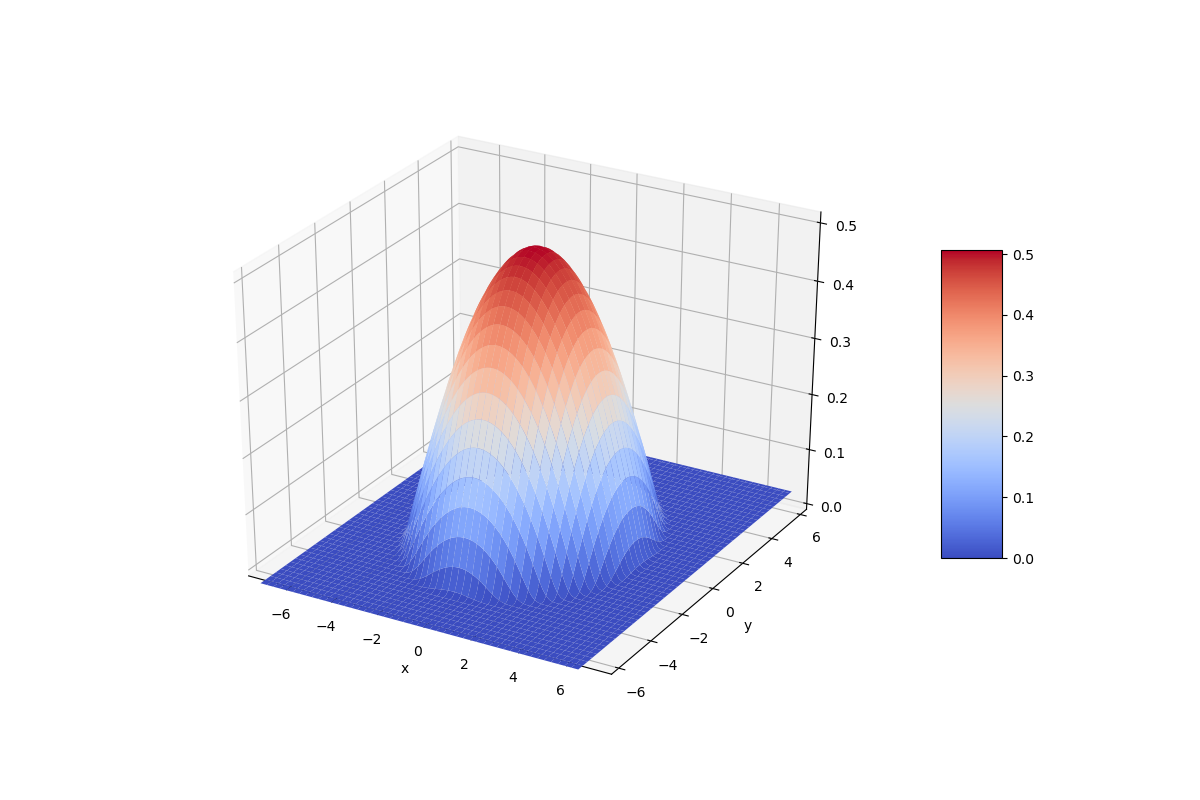}
         \caption{Learned solution slice}
         
     \end{subfigure}
     \hfill
    \begin{subfigure}[t]{0.3\textwidth}
         \centering
         \includegraphics[width=1.2\textwidth]{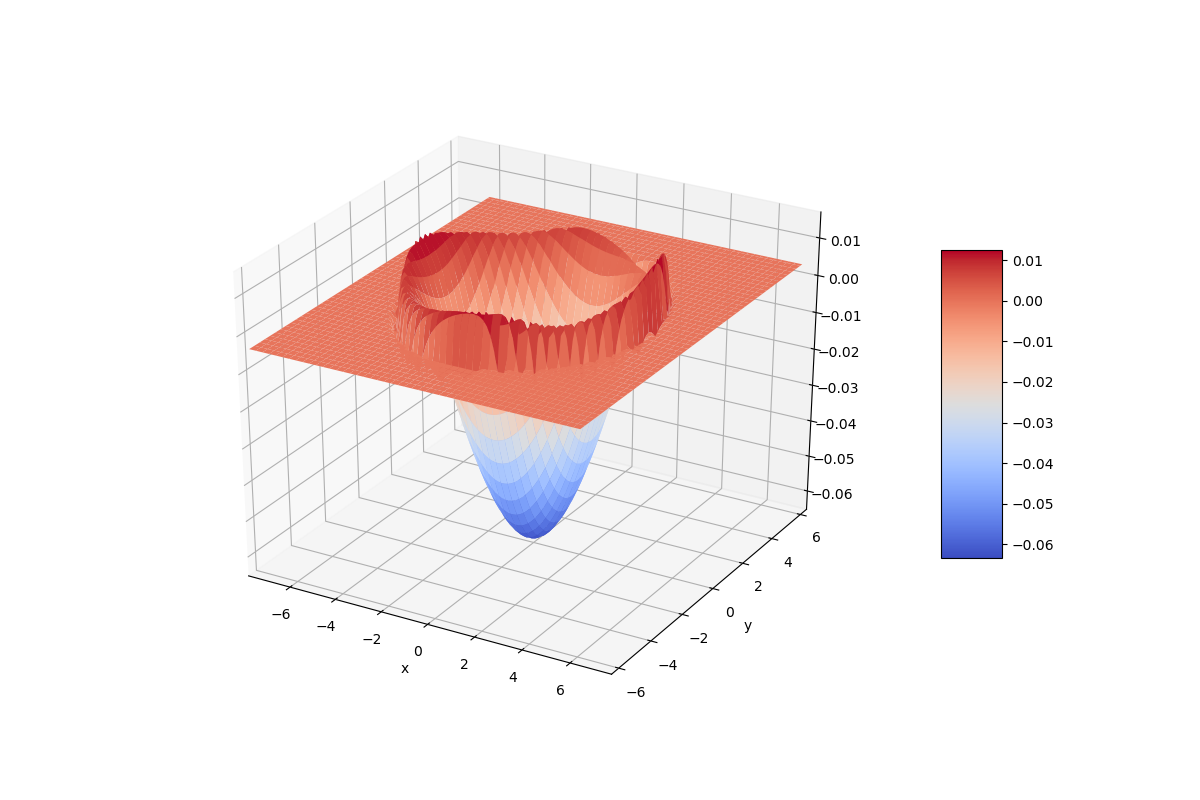}
         \caption{Learned solution error}
         
     \end{subfigure}
    \hfill
    \begin{subfigure}[t]{0.3\textwidth}
         \centering
         \includegraphics[width=1.2\textwidth]{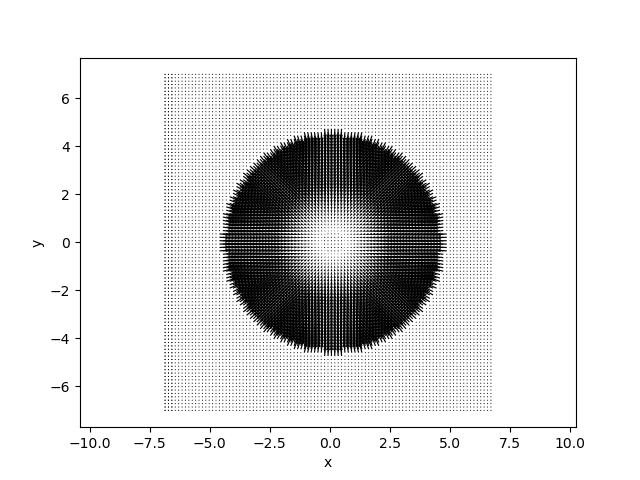}
         \caption{Barenblatt reference solution gradient}
         
     \end{subfigure}
     \hfill
    \begin{subfigure}[t]{0.3\textwidth}
         \centering
         \includegraphics[width=1.2\textwidth]{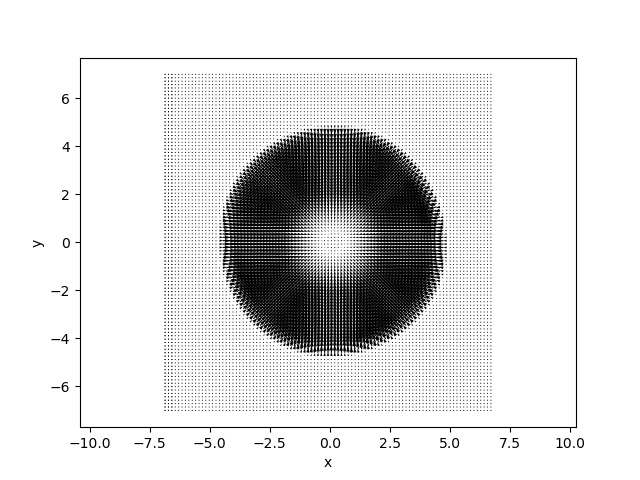}
         \caption{Learned solution gradient}
         
     \end{subfigure}
      \hfill
       \begin{subfigure}[t]{0.3\textwidth}
         \centering
         \includegraphics[width=1.2\textwidth]{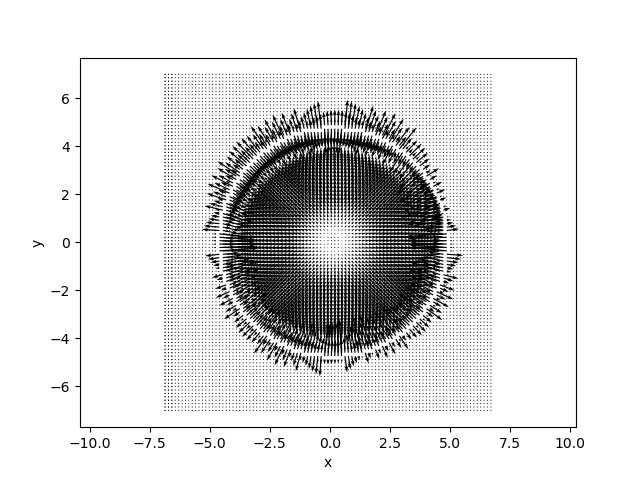}
         \caption{Learned solution gradient error}
         
     \end{subfigure}
        \caption{\textbf{15D, $L^1-$ PINN formulation  \eqref{PINN_full}:} Predicted solution slice  $u(0.5,x,y,1.0,\cdots, 1.0)$ for $\x\in \T = [-7,7]^{15}$, $t= 0.5$. }
\label{fig:PINN_l1}
\end{figure}

\begin{figure}
    \centering
    \begin{subfigure}[t]{0.47\textwidth}
             \centering
    \includegraphics[width = 1.2\textwidth]{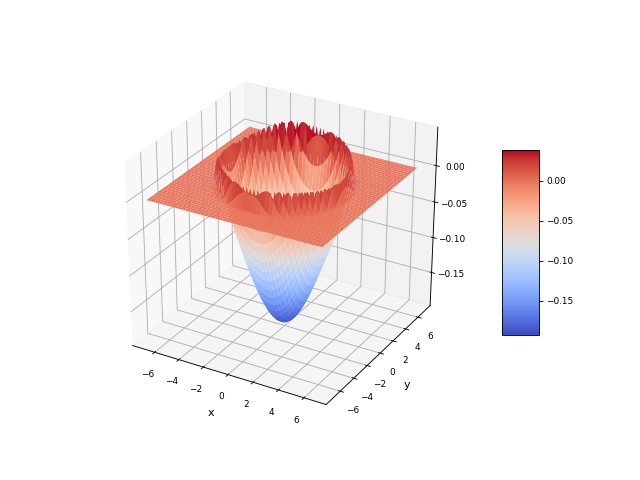}
               \caption{Learned $u_t$}
    \end{subfigure}
    \hfill
    \begin{subfigure}[t]{0.47\textwidth}
             \centering
    \includegraphics[width = 1.2\textwidth]{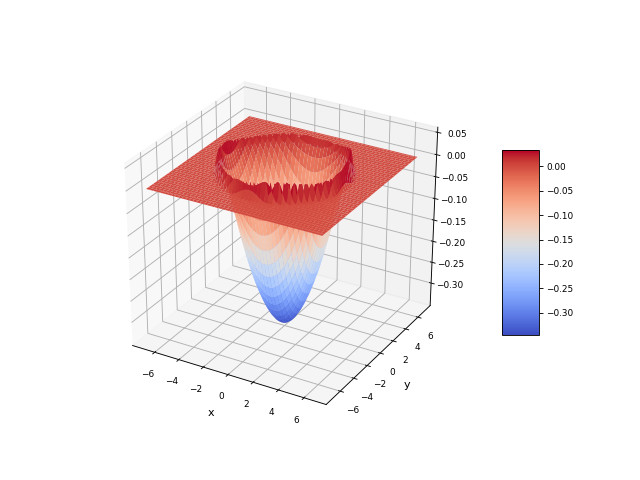}
               \caption{Learned $\displaystyle \frac{1}{2}\Delta u^2$}
    \end{subfigure}
    \caption{\textbf{15D, $L^1-$ PINN formulation \eqref{PINN_full}:} predicted partial derivatives.}
    \label{fig:PDE_l1}
\end{figure}

\begin{figure}
    \centering
    \begin{subfigure}[t]{0.47\textwidth}
             \centering
    \includegraphics[width = 1.2\textwidth]{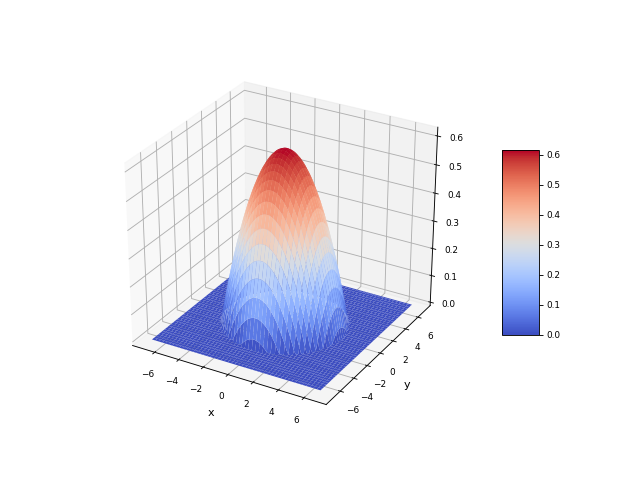}
               \caption{Exact $u_0$}
    \end{subfigure}
    \hfill
    \begin{subfigure}[t]{0.47\textwidth}
             \centering
    \includegraphics[width = 1.2\textwidth]{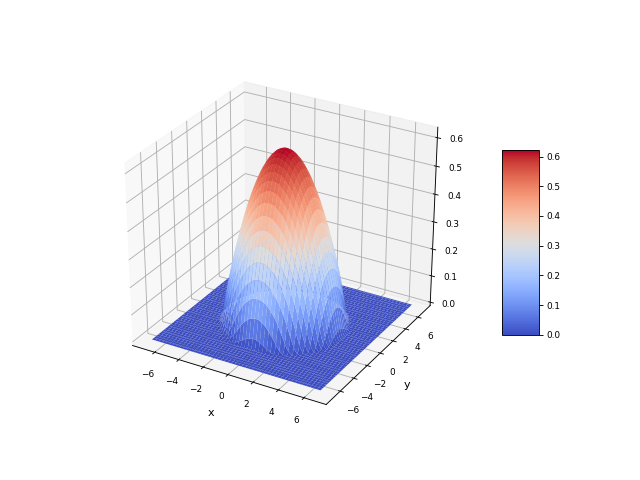}
               \caption{Learned initial value}
    \end{subfigure}
    
    \caption{\textbf{15D, $L^1-$ PINN formulation \eqref{PINN_full}}, predicted initial value $u(0,x,y,1.0,\cdots,1.0)$ for $\x \in 
    \T =[-7,7]^{15}$.}
    \label{fig:init_l1}
\end{figure}

\begin{figure}
    \centering
     \begin{subfigure}[t]{0.45\textwidth}
             \centering
    \includegraphics[width = 1.2\textwidth]{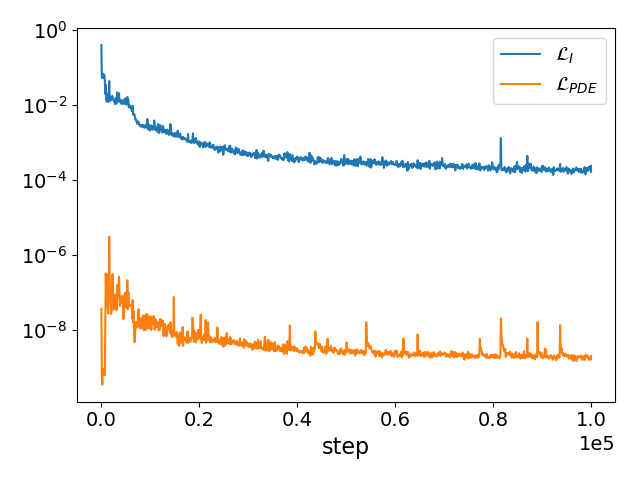}
               \caption{$L^2-$PINN}
    \end{subfigure}
    \hfill
     \begin{subfigure}[t]{0.45\textwidth}
             \centering
    \includegraphics[width = 1.2\textwidth]{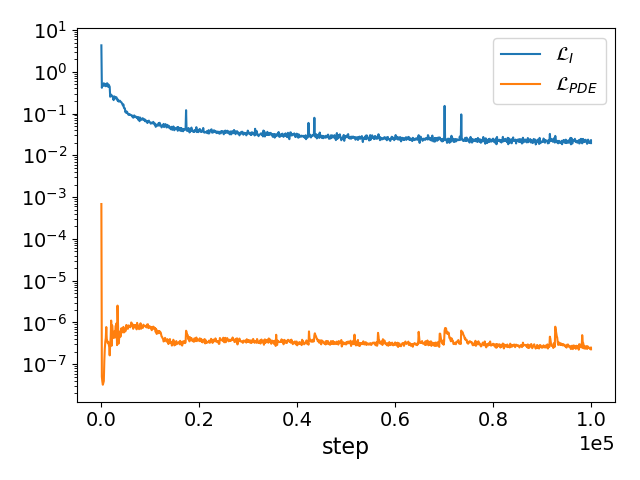}
    \caption{$L^1-$PINN}
    \end{subfigure}
    \caption{\textbf{15D:}  Training loss history by term.}
    \label{fig:training_history_PINN}
\end{figure}

From Table \ref{tab:PINN_error}, Figure \ref{fig:PINN_l2} and Figure \ref{fig:PINN_l1}, one can observe that the PINN formulation can indeed provide numerical solutions that closely approximate the exact ones even in high dimensions.
Not only is the neural network able to accurately approximate the function itself, but also the derivative of it.  This is essentially a result of successful imposition of the PDE. As can be observed from Figure \ref{fig:PDE_l2} and Figure \ref{fig:PDE_l1}, the learned $\partial_t u$ coincides with $\frac{1}{2} \Delta u^2$ which confirms the PDE has been successfully learned. The initial condition can also be softly enforced with term $\mathcal{L}_{I}$ as training proceeds (See Figure \ref{fig:init_l2} and Figure \ref{fig:init_l1} for an illustration). 
The training loss history of PINN is further presented in Figure \ref{fig:training_history_PINN} term by term. A training convergence can be observed from these plots, which further suggests a convergence to the exact solution ensured by \eqref{convergence_L1} and \eqref{convergence_L2}.

In addition, one can observe that learning the solution to QPME does not require the number of the trainable parameters of the solution ansatz to scale exponentially, which, in contrast to mesh-based solvers, is advantageous in dealing with high-d problems. 


\begin{table}[htbp]
    \hspace*{-1cm}
  \centering
    \begin{tabular}{|c|l|c|c|c|c|c|c|c|c|c|}
    \toprule
          \textbf{Dimension}&  & \multicolumn{1}{r|}{\textbf{1}} & \multicolumn{1}{r|}{\textbf{2}} & \multicolumn{1}{r|}{\textbf{3}} & \multicolumn{1}{r|}{\textbf{4}} & \multicolumn{1}{r|}{\textbf{5}} & \multicolumn{1}{r|}{\textbf{10}} & \multicolumn{1}{r|}{\textbf{15}} & \multicolumn{1}{r|}{\textbf{20}} & \multicolumn{1}{r|}{\textbf{50}} \\
    \midrule
    \multirow{3}[6]{3cm}{\textbf{Relative Error(\%)} for $L^2$-\textbf{PINN}} & \boldmath{}\textbf{$L^2$}\unboldmath{} & \multicolumn{1}{r|}{\textbf{0.21}} & \multicolumn{1}{r|}{\textbf{0.65}} & \multicolumn{1}{r|}{\textbf{0.61}} & \multicolumn{1}{r|}{\textbf{0.55}} & \multicolumn{1}{r|}{\textbf{1.1}} & \multicolumn{1}{r|}{\textbf{0.86}} & \multicolumn{1}{r|}{\textbf{1.72}} & \multicolumn{1}{r|}{\textbf{5.50}} & \multicolumn{1}{r|}{\textbf{16.03}} \\
\cmidrule{2-11}          & \boldmath{}\textbf{$L^1$}\unboldmath{} & \multicolumn{1}{r|}{\textbf{0.14}} & \multicolumn{1}{r|}{\textbf{0.4}} & \multicolumn{1}{r|}{\textbf{0.39}} & \multicolumn{1}{r|}{\textbf{0.43}} & \multicolumn{1}{r|}{\textbf{0.94}} & \multicolumn{1}{r|}{\textbf{0.71}} & \multicolumn{1}{r|}{\textbf{1.64}} & \multicolumn{1}{r|}{\textbf{5.12}} & \multicolumn{1}{r|}{\textbf{15.26}} \\
\cmidrule{2-11}          & \boldmath{}\textbf{$H^1$}\unboldmath{} & \multicolumn{1}{r|}{\textbf{4.28}} & \multicolumn{1}{r|}{\textbf{8.25}} & \multicolumn{1}{r|}{\textbf{7.93}} & \multicolumn{1}{r|}{\textbf{7.45}} & \multicolumn{1}{r|}{\textbf{8.64}} & \multicolumn{1}{r|}{\textbf{8.09}} & \multicolumn{1}{r|}{\textbf{9.88}} & \multicolumn{1}{r|}{\textbf{10.87}} & \multicolumn{1}{r|}{\textbf{28.5}} \\
    \midrule
    
 \multirow{3}[6]{3cm}{\textbf{Relative Error(\%) for $L^1$-PINN}} & \boldmath{}\textbf{$L^2$}\unboldmath{} & \multicolumn{1}{r|}{\textbf{0.45}} & \multicolumn{1}{r|}{\textbf{0.78}} & \multicolumn{1}{r|}{\textbf{0.95}} & \multicolumn{1}{r|}{\textbf{1.27}} & \multicolumn{1}{r|}{\textbf{2.07}} & \multicolumn{1}{r|}{\textbf{4.86}} & \multicolumn{1}{r|}{\textbf{10.46}} & \multicolumn{1}{r|}{\textbf{ 9.73}} & \multicolumn{1}{r|}{\textbf{10.76}} \\
\cmidrule{2-11}          & \boldmath{}\textbf{$L^1$}\unboldmath{} & \multicolumn{1}{r|}{\textbf{0.23}} & \multicolumn{1}{r|}{\textbf{0.5}} & \multicolumn{1}{r|}{\textbf{0.69}} & \multicolumn{1}{r|}{\textbf{1.1}} & \multicolumn{1}{r|}{\textbf{1.91}} & \multicolumn{1}{r|}{\textbf{4.65}} & \multicolumn{1}{r|}{\textbf{8.91}} & \multicolumn{1}{r|}{\textbf{8.37}} & \multicolumn{1}{r|}{\textbf{10.47}} \\
\cmidrule{2-11}          & \boldmath{}\textbf{$H^1$}\unboldmath{} & \multicolumn{1}{r|}{\textbf{5.11}} & \multicolumn{1}{r|}{\textbf{8.98}} & \multicolumn{1}{r|}{\textbf{9.38}} & \multicolumn{1}{r|}{\textbf{9.51}} & \multicolumn{1}{r|}{\textbf{10.48}} & \multicolumn{1}{r|}{\textbf{16.57}} & \multicolumn{1}{r|}{\textbf{16.45}} & \multicolumn{1}{r|}{\textbf{13.13}} & \multicolumn{1}{r|}{\textbf{21.55}} \\
\midrule
    \multirow{2}[4]{*}{Formulation Weight} & $\nu$ & \multicolumn{5}{c|}{$10^3$}             & \multicolumn{4}{c|}{1} \\
\cmidrule{2-11}          
& $\kappa$ & \multicolumn{5}{c|}{1}                & \multicolumn{4}{c|}{$10^3$} \\
    \midrule
    \multirow{2}[4]{*}{NN Architectire} & \# trainable & \multicolumn{1}{r|}{41001} & \multicolumn{1}{r|}{41201} & \multicolumn{1}{r|}{41401} & \multicolumn{1}{r|}{41601} & \multicolumn{1}{r|}{41801} & \multicolumn{1}{r|}{42801} & \multicolumn{1}{r|}{43801} & \multicolumn{1}{r|}{169601} & \multicolumn{1}{r|}{181601} \\
\cmidrule{2-11}          & Width/Depth & \multicolumn{7}{c|}{200/2}                            & \multicolumn{2}{c|}{400/2} \\
    \midrule
    \multirow{2}[4]{*}{Data Sampling} & $\theta_0$ & \multicolumn{9}{c|}{0.3} \\
\cmidrule{2-11}          & $\theta_1$ & \multicolumn{7}{c|}{0.3}                              & \multicolumn{2}{c|}{0.2} \\
    \midrule
    \multirow{2}[4]{*}{Training} &   Steps & \multicolumn{7}{c|}{$10^5$}     & \multicolumn{2}{c|}{$2\times10^5$}   \\
\cmidrule{2-11}          & Learning Rate & \multicolumn{9}{c|}{$10^{-3}$} \\
    \bottomrule
    \end{tabular}%
     \caption{\textbf{PINN formulation \eqref{PINN}: (hard Dirichlet B.C.+soft I.C.)} Relative error comparison for various dimensions .}
  \label{tab:PINN_error}%
\end{table}%

However, from Table \ref{tab:PINN_error}, one can observe that the generalization error of the neural network based solutions of high-dimensional cases are larger compared to that of low-dimensional cases. This could be a result of using larger neural network to approximate a more complicated solution in high-dimensional cases. Numerical experiments show that neural networks with width $200$ is no longer sufficient in approximating solutions to QPME with dimension larger than $15$, thus a larger network was adopted. Such network naturally requires more training and data for convergence. Since the training steps (data) was not quadrupled as the number of trainable does, which could have contributed to a larger approximation error. In addition, whether assigning quadrupled training steps and training data could bring significant improvement on the approximation accuracy is also questionable as the optimization of $\theta_{u}$ is highly nonconvex, which means one has to accept a significant and unavoidable uncertainty about optimization success with SGD or its variants.

\FloatBarrier

\subsubsection{\texorpdfstring{$\phi$}{} formulation}
In this section, we  consider the $\phi$ formulation \eqref{full_phi} to solve QPME \eqref{numerical_baren} in both $L^1$ and $L^2$ norm. The homogeneous Dirichlet boundary condition is enforced as a hard constraint following \eqref{phi_condition}. The initial condition can also be enforced softly with term $\mathcal{L}_{I}$ similar to the PINN formulation. The specific algorithmic settings are further presented in Table \ref{tab:phi_error} along with the relative errors computed for the trained solution slice $u(0.5,x,y,1.0,\cdots,1.0; \theta_u^*)$ at time $t= 0.5$ comparing with the exact solution \eqref{exact}.

\begin{figure}[htbp]
     \centering
     \begin{subfigure}[t]{0.3\textwidth}
         \centering
         \includegraphics[width=1.2\textwidth]{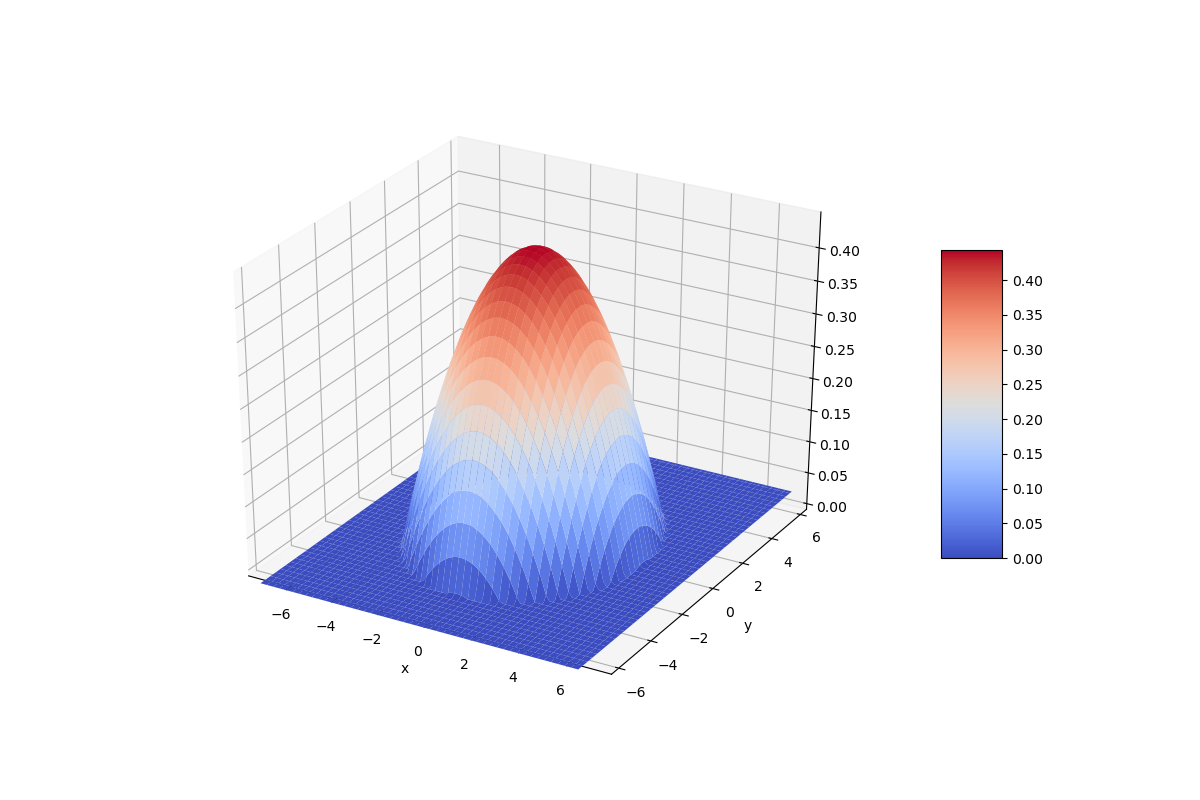}
         \caption{Barenblatt reference solution  }
         
     \end{subfigure}
    \hfill
     \begin{subfigure}[t]{0.3\textwidth}
         \centering
         \includegraphics[width=1.2\textwidth]{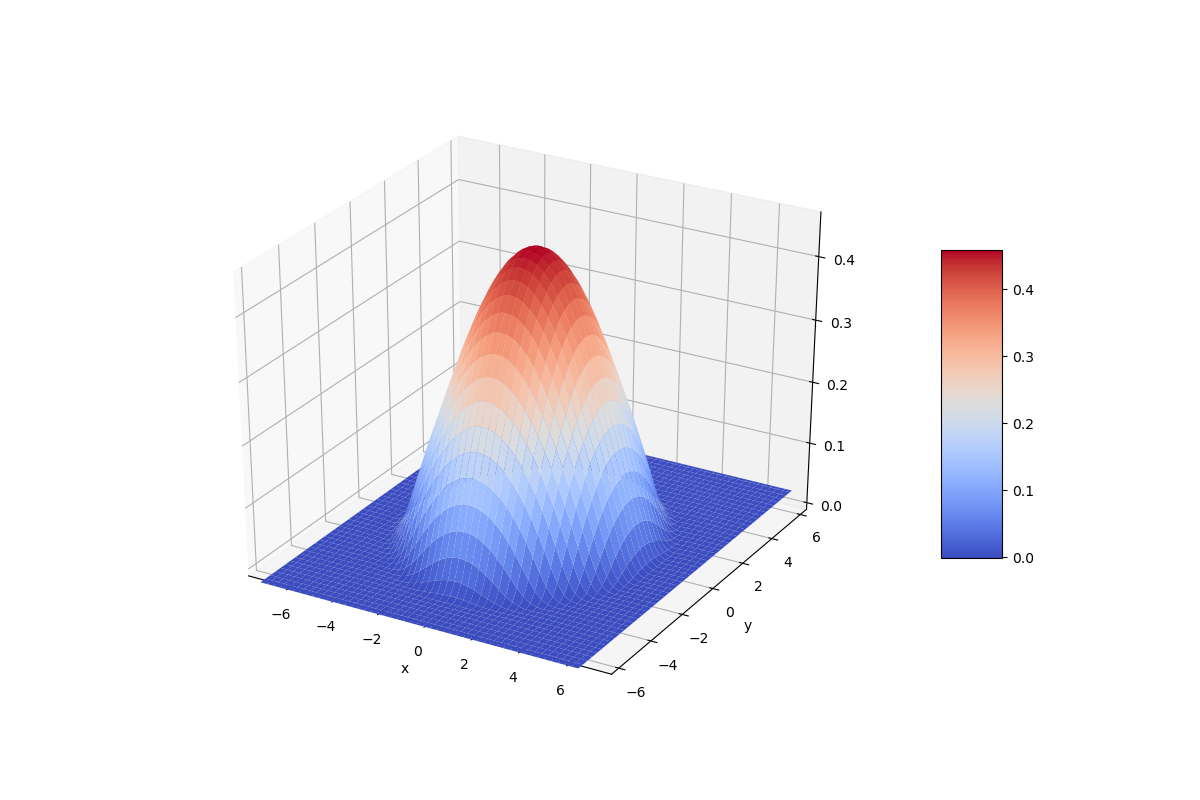}
         \caption{Learned solution slice}
         
     \end{subfigure}
     \hfill
    \begin{subfigure}[t]{0.3\textwidth}
         \centering
         \includegraphics[width=1.2\textwidth]{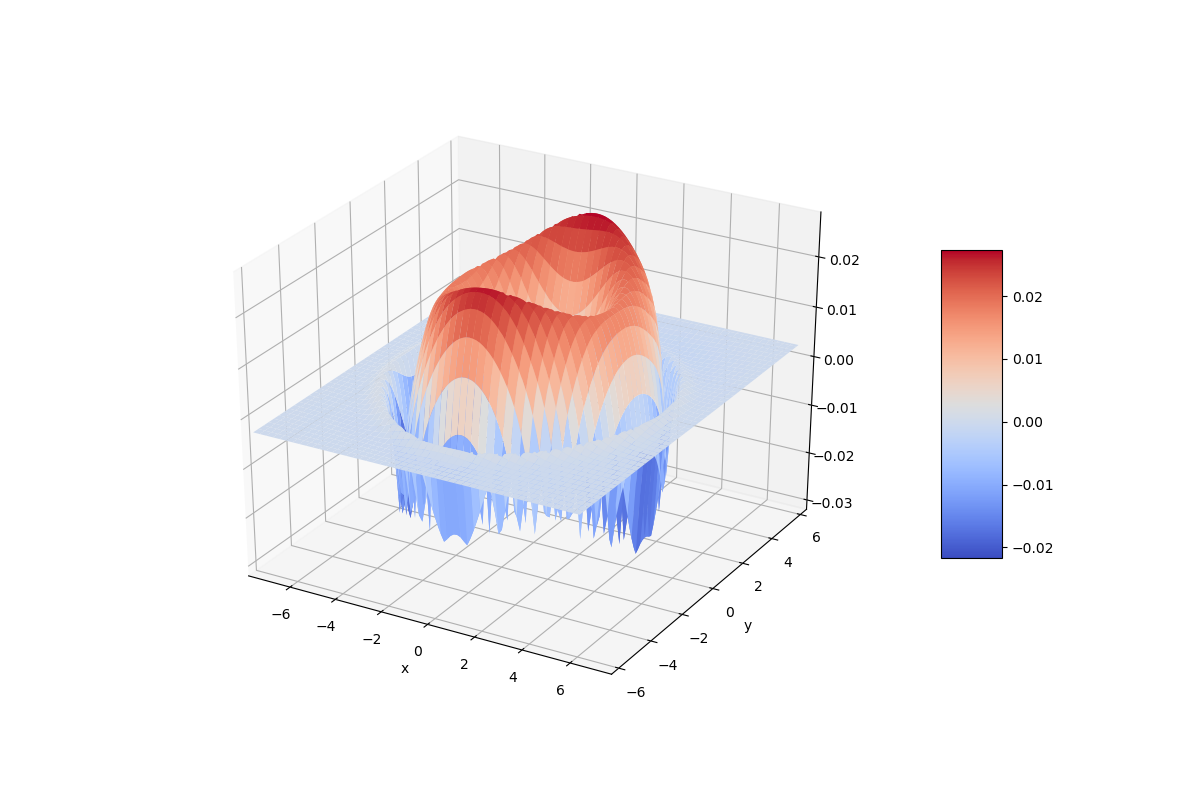}
         \caption{Learned solution error}
         
     \end{subfigure}
    \hfill
    \begin{subfigure}[t]{0.3\textwidth}
         \centering
         \includegraphics[width=1.2\textwidth]{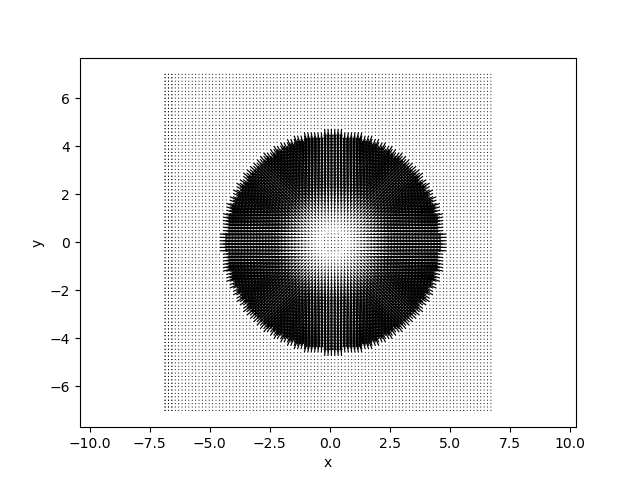}
         \caption{Barenblatt reference solution gradient}
         
     \end{subfigure}
     \hfill
    \begin{subfigure}[t]{0.3\textwidth}
         \centering
         \includegraphics[width=1.2\textwidth]{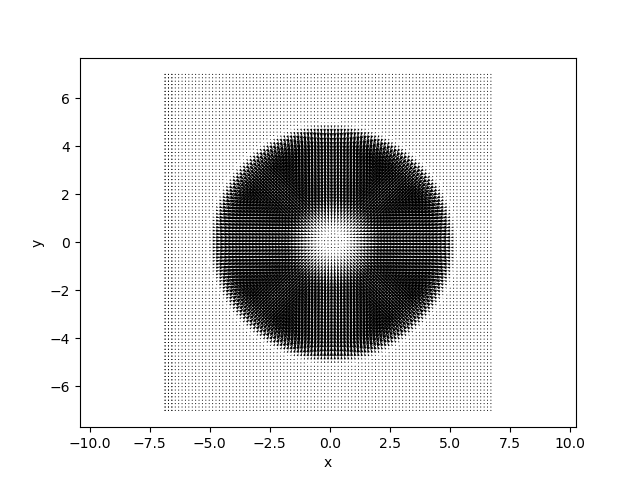}
         \caption{Learned solution gradient}
         
     \end{subfigure}
      \hfill
       \begin{subfigure}[t]{0.3\textwidth}
         \centering
         \includegraphics[width=1.2\textwidth]{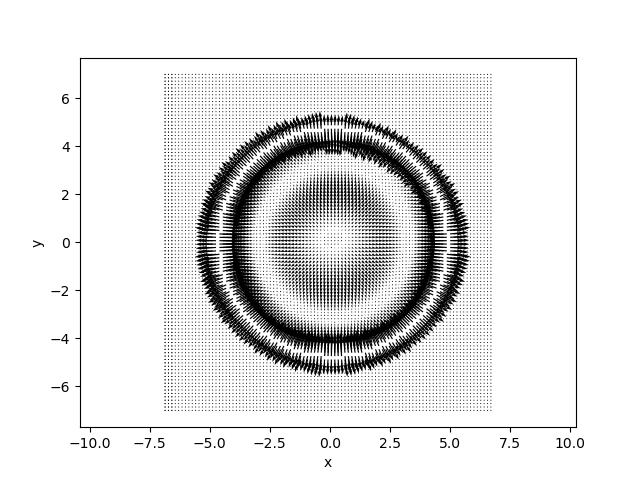}
         \caption{Learned solution gradient error}
         
     \end{subfigure}
        \caption{\textbf{15D, $L^2- \phi$  formulation  \eqref{phiformulation}:} Predicted solution slice  $u(0.5,x,y,1.0,\cdots, 1.0)$ for $\x\in \T = [-7,7]^{15}$, $t= 0.5$. }
\label{fig:phi_l2}
\end{figure}

\begin{figure}[htbp]
     \centering
     \begin{subfigure}[t]{0.3\textwidth}
         \centering
         \includegraphics[width=1.2\textwidth]{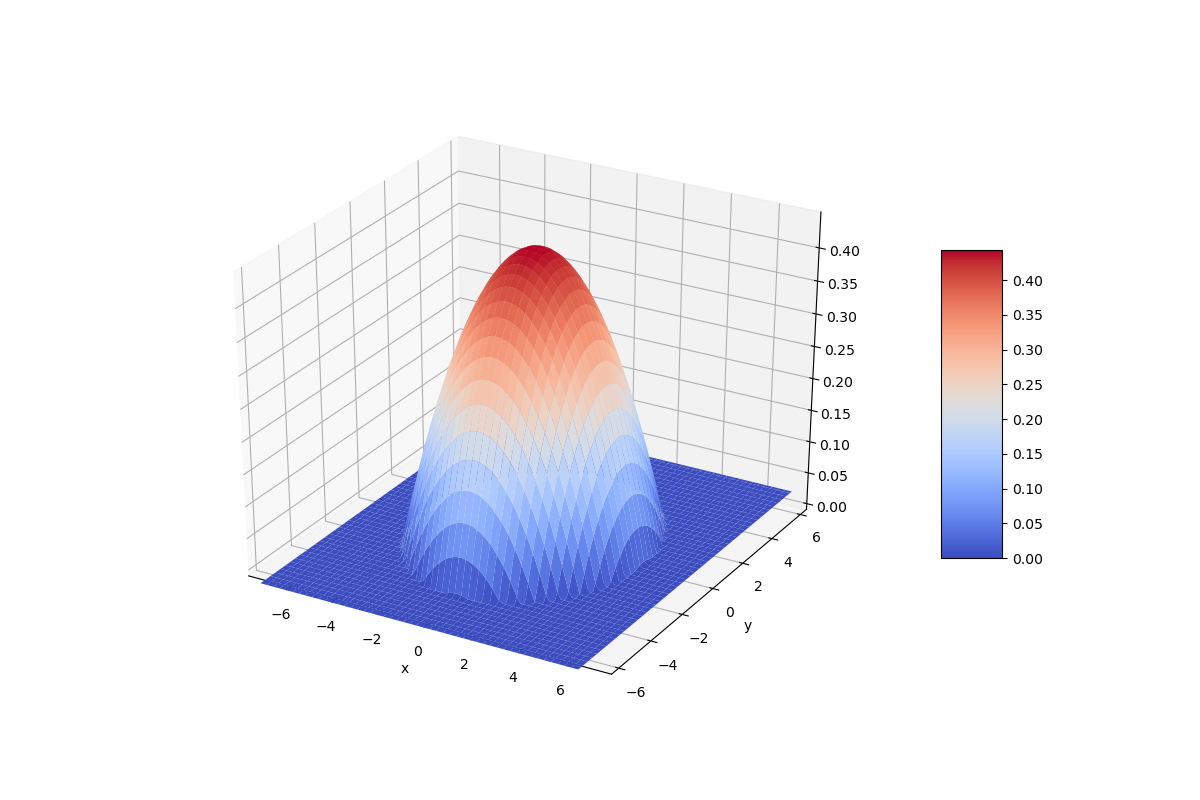}
         \caption{Barenblatt reference solution}
         
     \end{subfigure}
    \hfill
     \begin{subfigure}[t]{0.3\textwidth}
         \centering
         \includegraphics[width=1.2\textwidth]{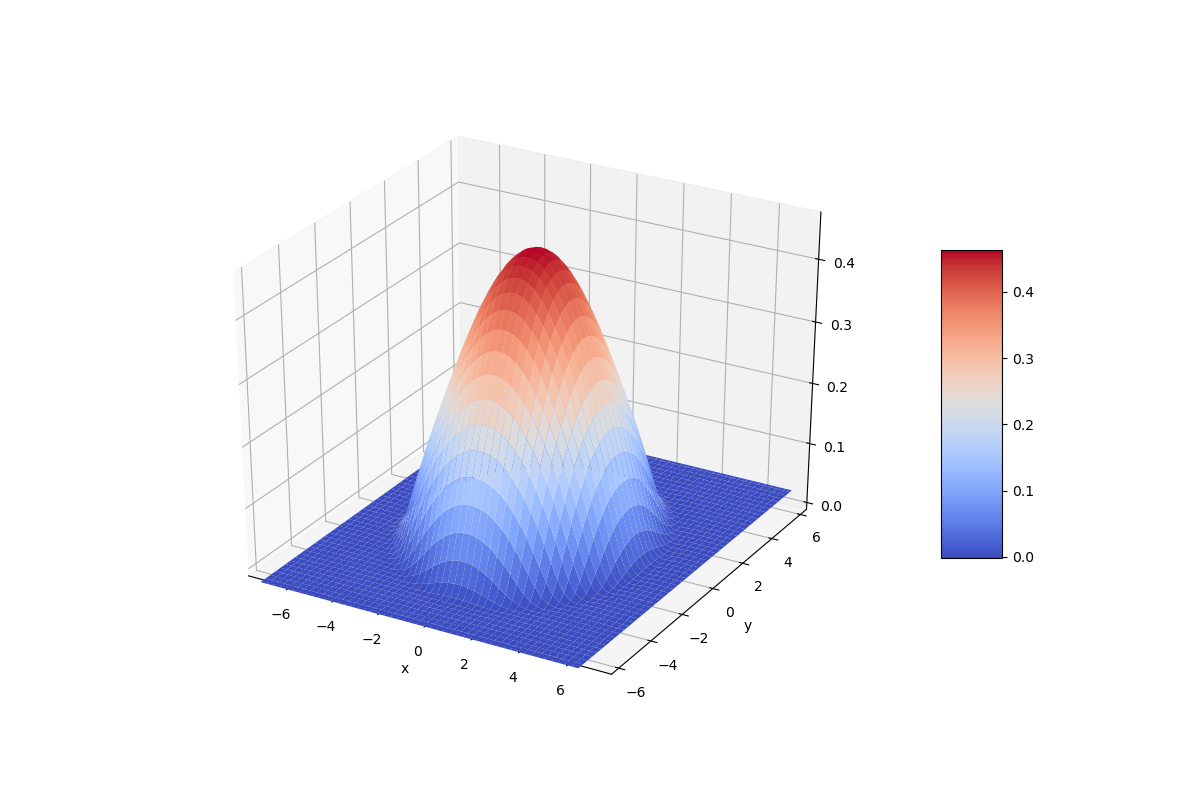}
         \caption{Learned solution slice}
     \end{subfigure}
     \hfill
    \begin{subfigure}[t]{0.3\textwidth}
         \centering
         \includegraphics[width=1.2\textwidth]{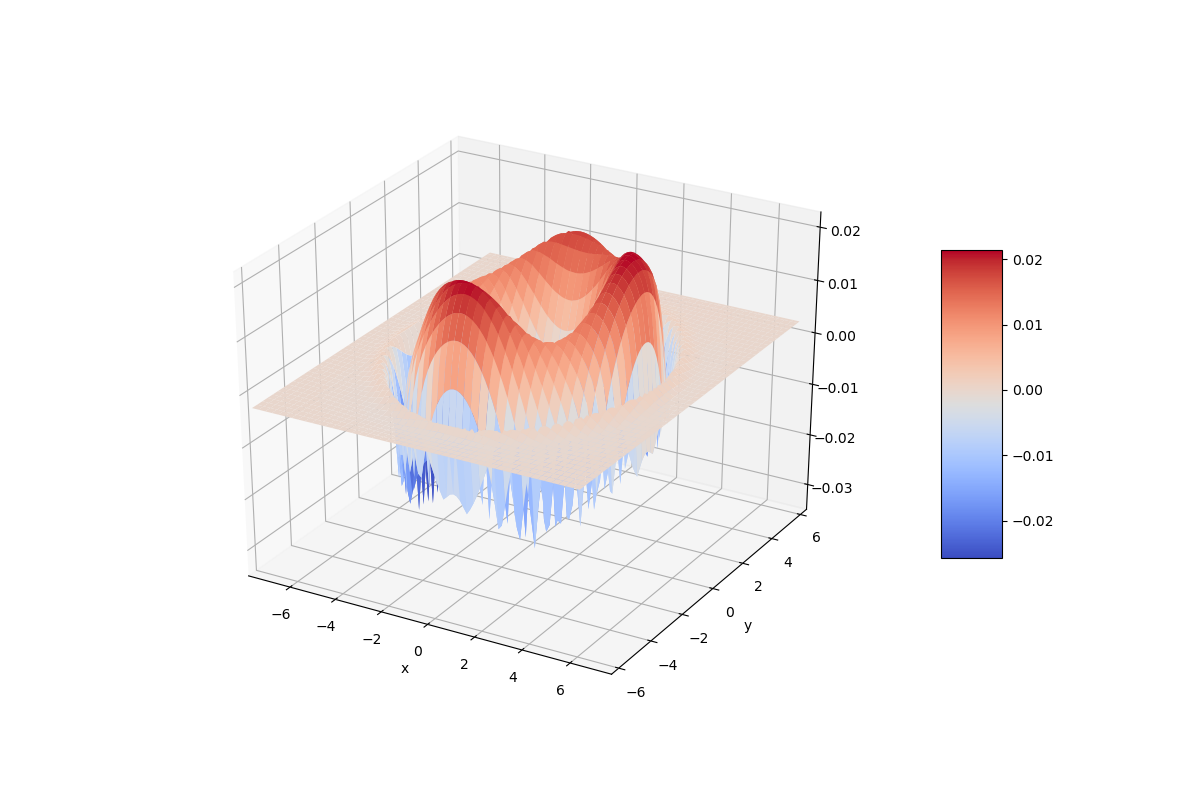}
         \caption{Learned solution error}
     \end{subfigure}
    \hfill
    \begin{subfigure}[t]{0.3\textwidth}
         \centering
         \includegraphics[width=1.2\textwidth]{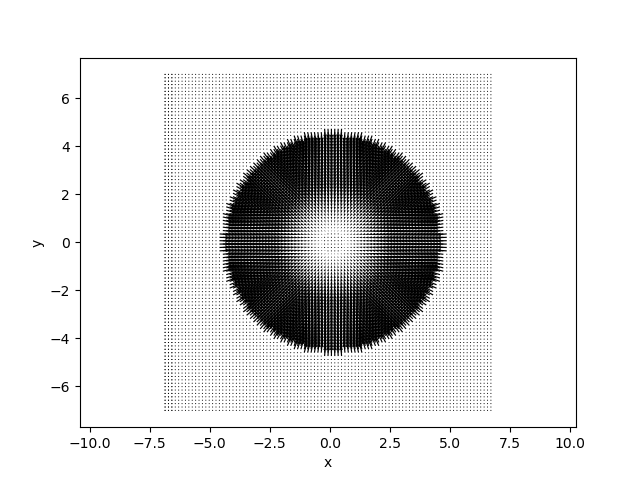}
         \caption{Barenblatt reference solution gradient}
     \end{subfigure}
     \hfill
    \begin{subfigure}[t]{0.3\textwidth}
         \centering
         \includegraphics[width=1.2\textwidth]{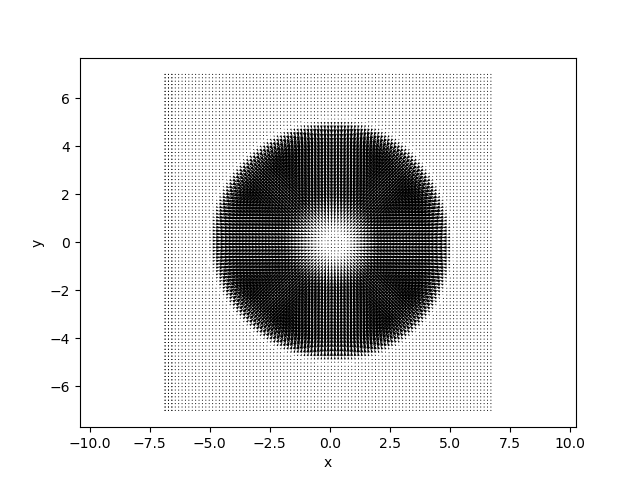}
         \caption{Learned solution gradient}
         
     \end{subfigure}
      \hfill
       \begin{subfigure}[t]{0.3\textwidth}
         \centering
         \includegraphics[width=1.2\textwidth]{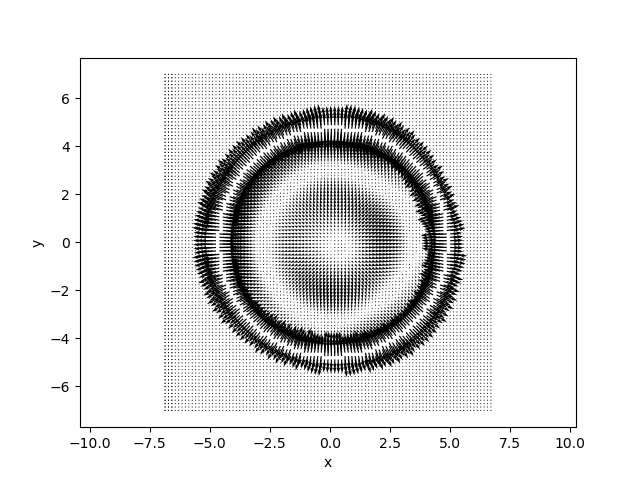}
         \caption{Learned solution gradient error}
         
     \end{subfigure}
        \caption{\textbf{15D,  $L^1-\phi$ formulation \eqref{phiformulation}:} Predicted solution slice  $u_{\phi}(0.5,x,y,1.0,\cdots, 1.0)$ with  for $\x\in \T = [-7,7]^{15}$, $t= 0.5$. }
\label{fig:phi_l1}
\end{figure}

\begin{figure}
    \centering
    \begin{subfigure}[t]{0.45\textwidth}
        \includegraphics[width = 1.2\textwidth]{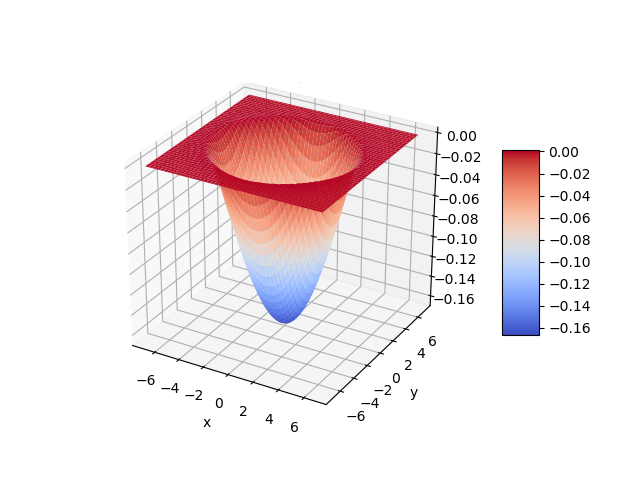}
        \caption{Learned through $L^1$-$\phi$ formulation}
    \end{subfigure}
    \hfill
     \begin{subfigure}[t]{0.45\textwidth}
        \includegraphics[width = 1.2\textwidth]{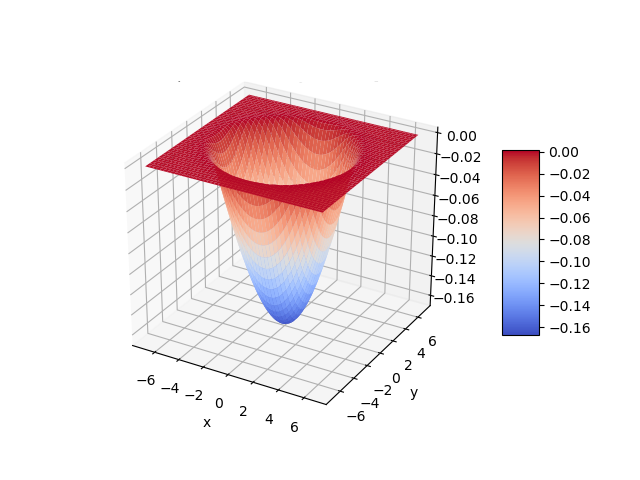}
        \caption{Learned through $L^2$-$\phi$ formulation}
    \end{subfigure}
    \caption{\textbf{15D: }Predicted $\phi(0.5, x,y,1.0,\cdots,1.0; \theta_{\phi}^*)$.}
    \label{fig:phi_pred}
\end{figure}

From Table \ref{tab:phi_error}, Figure \ref{fig:phi_l2} and Figure \ref{fig:phi_l1}, one can observe that the $\phi$ formulation can indeed provide numerical solutions that closely approximate the exact ones up till dimension $20$.
Not only is the neural network able to accurately approximate the function itself, but also the derivatives of it. The mismatch mainly concentrated near the region where the solution is not smooth (free boundary). 
The predicted minimizer $\phi$ to \eqref{full_phi} is also presented as in Figure \ref{fig:phi_pred}. 


Theoretically speaking, compared to PINN, $\phi$ formulation is advantageous as it can be applied to solve a wider range of QPMEs, whose solution are less regular or smooth. 
However, for the case being tested, we do encounter more challenges in the training process especially for the high-dimensional cases compared with that of PINN. 
One observation is that the generalization error of the testing solution slice gets larger as the dimension gets higher. This could attributes to the nature of the exact solution $U_2$ as its nonzero region only accounts for a tiny small portion of $\T$($\ll 1$\textperthousand) when $d$ is large.
That is to say, the zero function will be a pretty good approximation of $U_2$ already in both $L^1(\T)$ and $L^2 (\T)$ sense. The training can thus be easily trapped in a local minimum $u_{\phi} =0$ which can be reached by a neural network $\phi = 0$. In addition, the generalization error reported measures the error for a solution slice projected onto a two-dimensional space in stead of in $\T$ to ease the computation and visualization, which can be an uncomprehensive measurement of the error. Moreover, the selected slice is a slice whose values are dominated by nonzero ones, which can also be a unfair representative of the entire solution to quantify the relative error.

The reason that PINN formulation seems to suffer less from this effects is probably the application of efficient sampling. Since the correction term $c(\x)$ is {{not}} applied for any terms in the loss functional of PINN, meaning a very large weight was put on the region where the solution is nonzero when evaluating $\mathcal{L}_{\text{PINN}}$ , which could have helped the solution ansatz to escape from the local minimum. However, by the nature of $\phi$ formulation, the $c(\x)$ can not be omitted. Otherwise, the target functional will be changed entirely.

Furthermore, while we are able to identify the desired solutions in many cases, theoretically, one can not guarantee meaningful solutions to QPME form the training of $\phi$ formulation. In fact, both the condition $1-\Delta \phi$ and $u_{\phi}\geq 0$ are not enforced at all in this formulation. Theses conditions can only be used post training to carry out a solution selection or as a criteria for early truncation of training.


Artificial choices of other algorithmic ingredients such as batch sizes, learning rates, $\theta_0$ and $\theta_1$ will also inevitable influence the optimization process providing limited computational resources.


\begin{table}[htbp]
  \centering
  \hspace*{-1.8cm}
    \begin{tabular}{|c|l|c|c|c|c|c|c|c|c|r|}
    \toprule
         \textbf{Dimension} &  & \multicolumn{1}{r|}{\textbf{1}} & \multicolumn{1}{r|}{\textbf{2}} & \multicolumn{1}{r|}{\textbf{3}} & \multicolumn{1}{r|}{\textbf{4}} & \multicolumn{1}{r|}{\textbf{5}} & \multicolumn{1}{r|}{\textbf{10}} & \multicolumn{1}{r|}{\textbf{15}} & \multicolumn{1}{r|}{\textbf{20}} & \textbf{50} \\
    \midrule
    \multirow{3}[6]{3.5cm}{\textbf{Relative Errors(\%) for $L^2-\phi$ Formulation}} & \boldmath{}\textbf{$L^2$}\unboldmath{} & \multicolumn{1}{r|}{\textbf{3.58}} & \multicolumn{1}{r|}{\textbf{4.95}} & \multicolumn{1}{r|}{\textbf{4.41}} & \multicolumn{1}{r|}{\textbf{9.77}} & \multicolumn{1}{r|}{\textbf{5.77}} & \multicolumn{1}{r|}{\textbf{3.82}} & \multicolumn{1}{r|}{\textbf{8.29}} & \multicolumn{1}{r|}{\textbf{14.45}} & \textbf{54.26} \\
\cmidrule{2-11}          & \boldmath{}\textbf{$L^1$}\unboldmath{} & \multicolumn{1}{r|}{\textbf{3.23}} & \multicolumn{1}{r|}{\textbf{5.87}} & \multicolumn{1}{r|}{\textbf{4.57}} & \multicolumn{1}{r|}{\textbf{9.98}} & \multicolumn{1}{r|}{\textbf{6.22}} & \multicolumn{1}{r|}{\textbf{3.97}} & \multicolumn{1}{r|}{\textbf{8.99}} & \multicolumn{1}{r|}{\textbf{15.56}} & \textbf{77.65} \\
\cmidrule{2-11}          & \boldmath{}\textbf{$H^1$}\unboldmath{} & \multicolumn{1}{r|}{\textbf{16.44}} & \multicolumn{1}{r|}{\textbf{18.84}} & \multicolumn{1}{r|}{\textbf{20.18}} & \multicolumn{1}{r|}{\textbf{27.29}} & \multicolumn{1}{r|}{\textbf{16.97}} & \multicolumn{1}{r|}{\textbf{14.98}} & \multicolumn{1}{r|}{\textbf{19.21}} & \multicolumn{1}{r|}{\textbf{24.15}} & \textbf{71} \\
\midrule
    \multirow{3}[6]{3.5cm}{\textbf{Relative Errors(\%) for $L^1-\phi$ Formulation}} & \boldmath{}\textbf{$L^2$}\unboldmath{} & \multicolumn{1}{r|}{\textbf{2.32}} & \multicolumn{1}{r|}{\textbf{5.01}} & \multicolumn{1}{r|}{\textbf{5.06}} & \multicolumn{1}{r|}{\textbf{8.8}} & \multicolumn{1}{r|}{\textbf{4.87}} & \multicolumn{1}{r|}{\textbf{3.62}} & \multicolumn{1}{r|}{\textbf{9.25}} & \multicolumn{1}{r|}{\textbf{26.67}} & \textbf{52.24} \\
\cmidrule{2-11}          & \boldmath{}\textbf{$L^1$}\unboldmath{} & \multicolumn{1}{r|}{\textbf{2.11}} & \multicolumn{1}{r|}{\textbf{6.16}} & \multicolumn{1}{r|}{\textbf{5.84}} & \multicolumn{1}{r|}{\textbf{9.21}} & \multicolumn{1}{r|}{\textbf{4.84}} & \multicolumn{1}{r|}{\textbf{3.52}} & \multicolumn{1}{r|}{\textbf{9.49}} & \multicolumn{1}{r|}{\textbf{28.05}} & \textbf{85.13} \\
\cmidrule{2-11}          & \boldmath{}\textbf{$H^1$}\unboldmath{} & \multicolumn{1}{r|}{\textbf{14.52}} & \multicolumn{1}{r|}{\textbf{25.25}} & \multicolumn{1}{r|}{\textbf{18.23}} & \multicolumn{1}{r|}{\textbf{25.73}} & \multicolumn{1}{r|}{\textbf{16.82}} & \multicolumn{1}{r|}{\textbf{13.73}} & \multicolumn{1}{r|}{\textbf{20.3}} & \multicolumn{1}{r|}{\textbf{44.74}} & \textbf{67.49} \\
\midrule
    \multirow{2}[4]{*}{Formulation Weights} & $\nu$ & \multicolumn{4}{c|}{$10^3$}     & \multicolumn{5}{c|}{1} \\
\cmidrule{2-11}          & $\kappa$ & \multicolumn{4}{c|}{1}        & \multicolumn{2}{c|}{$10^3$} &$10^4$ & \multicolumn{1}{r|}{$10^5$}& $10^3$ \\
    \midrule
    \multirow{2}[4]{*}{NN Architecture} & Width/Depth & \multicolumn{7}{c|}{200/2}                              & \multicolumn{2}{c|}{400/2} \\
\cmidrule{2-11}          & \# trainable & \multicolumn{1}{r|}{41001} & \multicolumn{1}{r|}{41201} & \multicolumn{1}{r|}{41401} & \multicolumn{1}{r|}{41601} & \multicolumn{1}{r|}{41801} & \multicolumn{1}{r|}{42801} & \multicolumn{1}{r|}{43801} & \multicolumn{1}{r|}{169601} & 181601 \\
    \midrule
    \multirow{2}[4]{*}{Data Sampling} & $\theta_0$ & \multicolumn{7}{c|}{0.3}                              & \multicolumn{1}{r|}{0.2} & 0.4 \\
\cmidrule{2-11}          & $\theta_1$ & \multicolumn{8}{c|}{0.3}                             & 0.4 \\
    \midrule
    \multirow{2}[4]{*}{Training} &   Steps & \multicolumn{7}{c|}{$10^5$}                   &
    \multicolumn{1}{r|}{$6\times 10^{5}$} & $2\times{10^{5}}$ \\
\cmidrule{2-11}          & Learning Rate & \multicolumn{8}{c|}{$10^{-3}$}                                    & $10^{-4}$\\
    \bottomrule
    \end{tabular}%
     \caption{\textbf{$\phi$ formulation \eqref{phiformulation} (hard Dirichlet B.C.+soft I.C.): } Relative error comparison for various dimensions.}
  \label{tab:phi_error}%
\end{table}%

\pagebreak

\begin{figure}[htbp]
    \centering
     \begin{subfigure}[t]{0.45\textwidth}
             \includegraphics[width =1.2\textwidth]{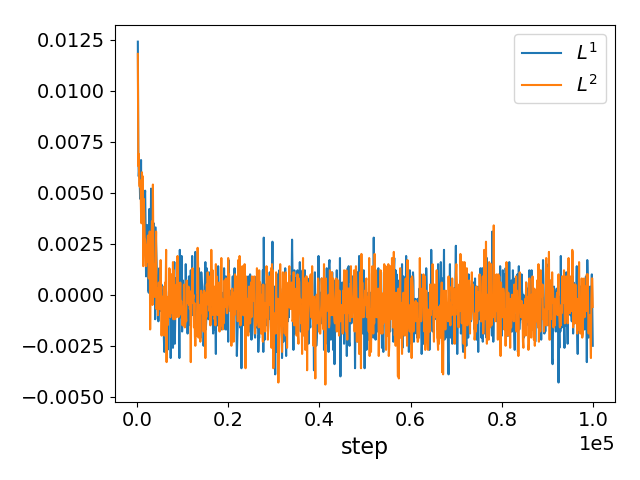}
             \caption{\textbf{5D}}
    \end{subfigure}
    \hfill
    \begin{subfigure}[t]{0.45\textwidth}
             \includegraphics[width =1.2\textwidth]{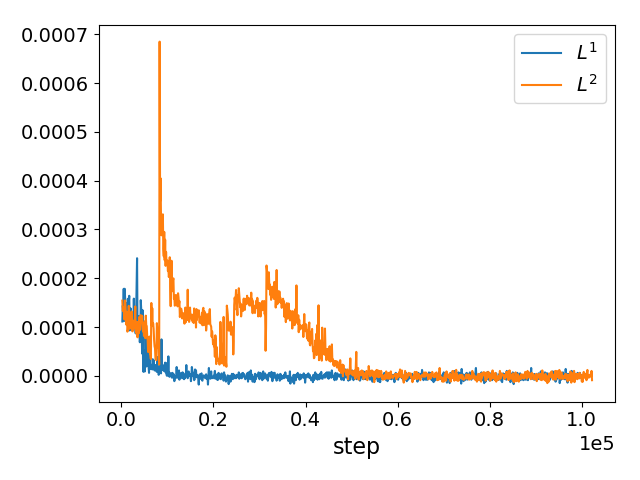}
             \caption{\textbf{10D}}
    \end{subfigure}
    \caption{Empirical $\int_Q U_2^2(t,\x) +\mathcal{L}_{\phi}\big(u_{\phi}(t,\x;\theta_{\phi} )\big)$ as training proceeds. }
    \label{fig:phi_usq_comp}
\end{figure}

We further observe that the optimization of $\mathcal{L}_{\phi}(u_{\phi}(t,\x;\theta_{\phi}))$ indeed converges to $-\int_{Q}U_2^2$ as training proceeds with $u_{\phi}$ being the parametrize solution ansatz as stated in \eqref{u_phi}.
This observation in fact confirms the theoretical result \eqref{form_equivalency} derived in \cite{brenier2020examples}.
In Figure \ref{fig:phi_usq_comp}, we specifically use the batch of training data at each training step to evaluate empirically the value of  $-\int_Q U_2^2$ for the exact solution $U_2(t,\x)$ defined as in \eqref{exact}, and further compare it with the empirical loss $\mathcal{L}_{\phi}(u_{\phi})$ based on the neural network solution $u_{\phi}$ at that time. As one can observe, the difference of the two values gradually reduces as the training continues, which verifies the training effectiveness of this formulation.

\FloatBarrier

\subsection{\texorpdfstring{$q$}{}-\texorpdfstring{$\sigma$}{} formulation}
Since $q-\sigma$ formulation \eqref{qsigma_full} is developed based on the $\phi$ formulation. The training thus also suffers from challenges met in the training of $\phi$ formulation, i.e., the training can be easily trapped in a local minimum, i.e. $u_{q,\sigma} = 0$ for high-dimensional cases. Additionally, the partial derivatives of $\phi$ is separated into two independent functions $q$ and $\sigma$, whose correlation is only enforced softly with the loss term $\mathcal{L}_{\partial_t \sigma, \Delta q}$, which can pose more challenges to the optimization of the target functional. Additional hyper-parameter $\gamma$ is also introduced to adjust the weight of $\mathcal{L}_{\partial_t \sigma, \Delta q}$ whose optimal choice is again obscure. Due to such reasons, only results for dimension $1$ to $10$ are reported as no reasonable results for higher-dimensions were obtained in the scope of experiments that were carried out.

\begin{figure}[htbp]
     \centering
     \begin{subfigure}[t]{0.3\textwidth}
         \centering
         \includegraphics[width=1.2\textwidth]{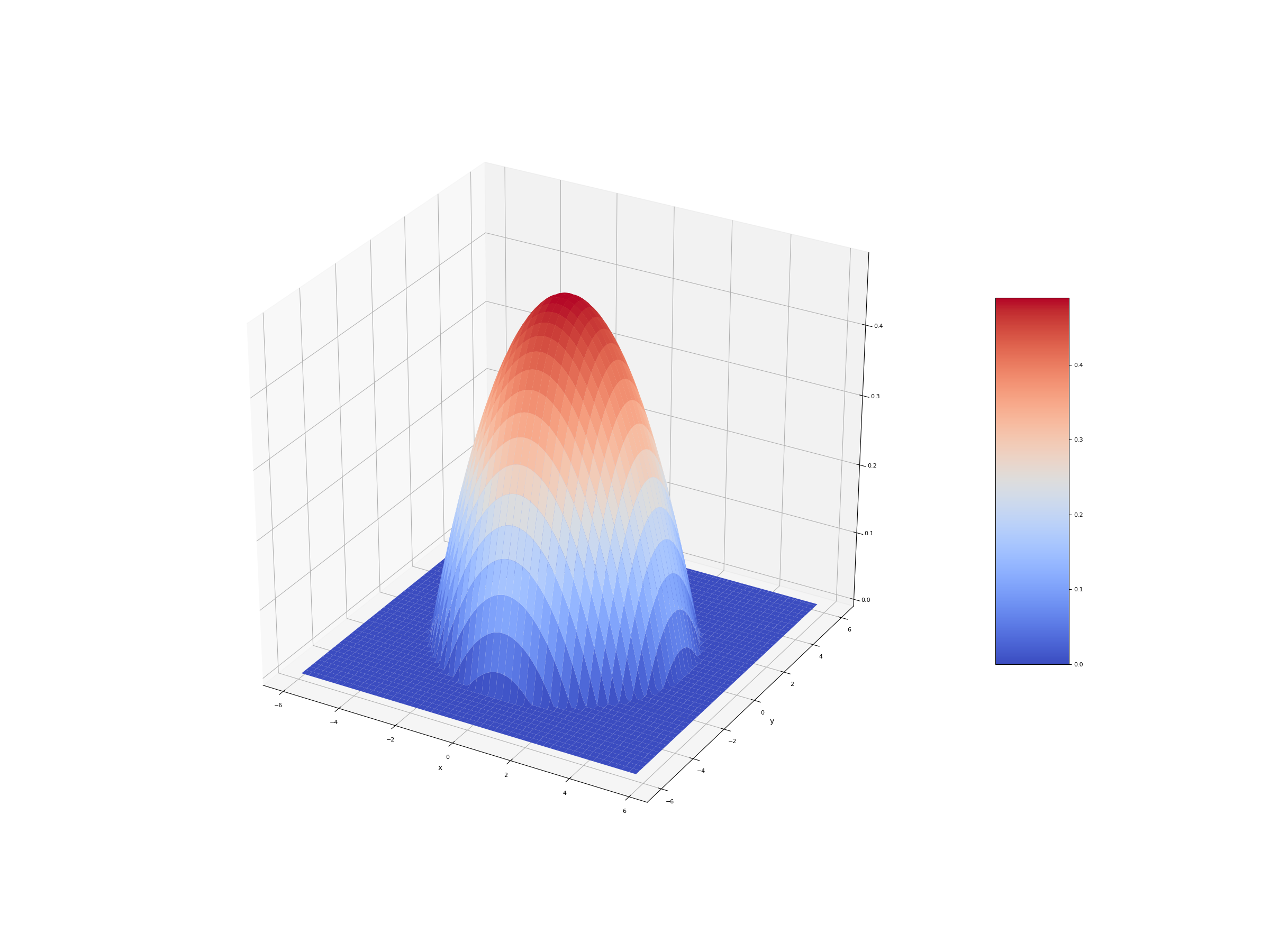}
         \caption{Barenblatt reference solution  }
         
     \end{subfigure}
    \hfill
     \begin{subfigure}[t]{0.3\textwidth}
         \centering
         \includegraphics[width=1.2\textwidth]{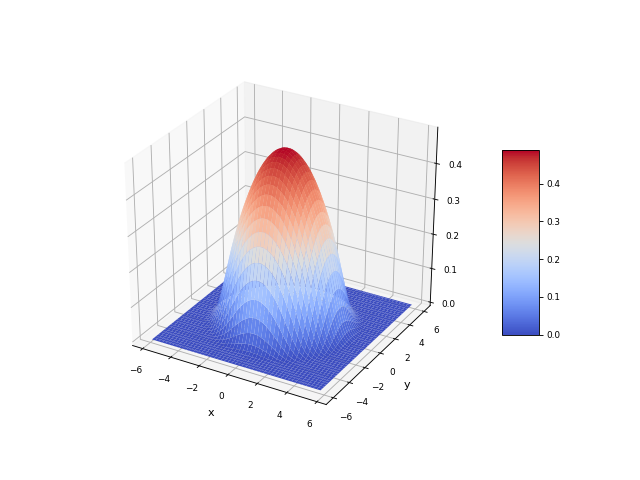}
         \caption{Learned solution slice}
         
     \end{subfigure}
     \hfill
    \begin{subfigure}[t]{0.3\textwidth}
         \centering
\includegraphics[width=1.2\textwidth]{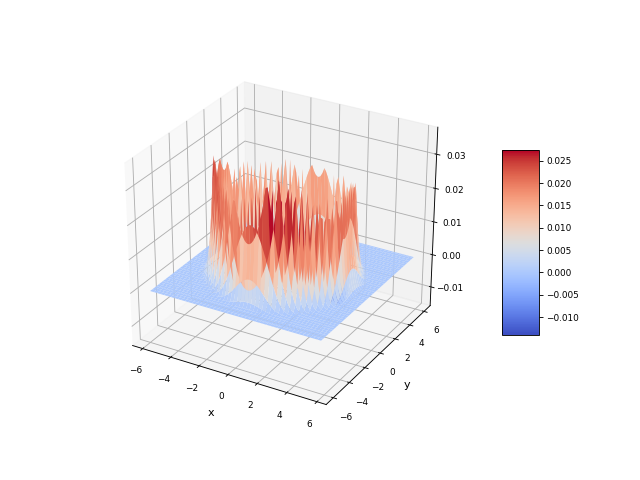}
         \caption{Learned solution error}
         
     \end{subfigure}
        \caption{\textbf{10D,  $L^2-q-\sigma$ formulation \eqref{qsigma_full}:} Predicted solution slice  $u(0.5,x,y,1.0,\cdots, 1.0)$ for $\x\in \T = [-6,6]^{10}$, $t= 0.5$. }
\label{fig:qsigma_l2}
\end{figure}

\begin{figure}
    \centering
    \begin{subfigure}[t]{0.47\textwidth}
      \includegraphics[width = 1.2\textwidth]{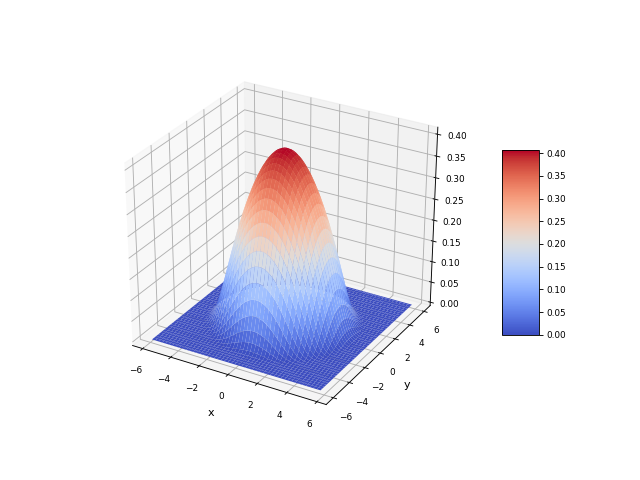}
    \caption{Learned $q$}
    \end{subfigure}
    \hfill
    \begin{subfigure}[t]{0.47\textwidth}
    \centering
    \includegraphics[width = 1.2\textwidth]{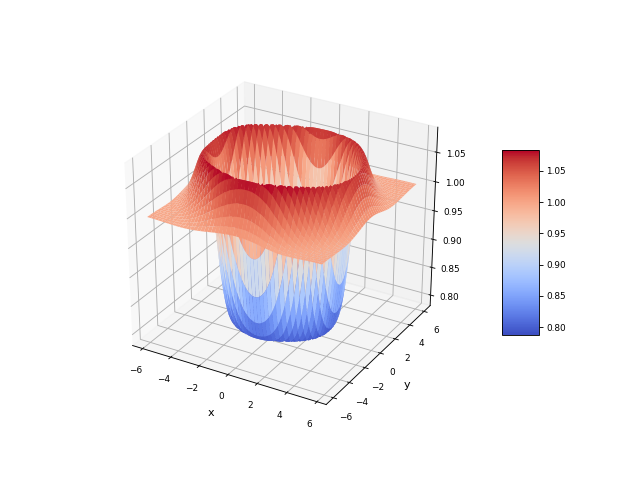}
    \caption{Learned $\displaystyle \sigma$}
    \end{subfigure}
    \begin{subfigure}[t]{0.47\textwidth}
             \centering
    \includegraphics[width = 1.2\textwidth]{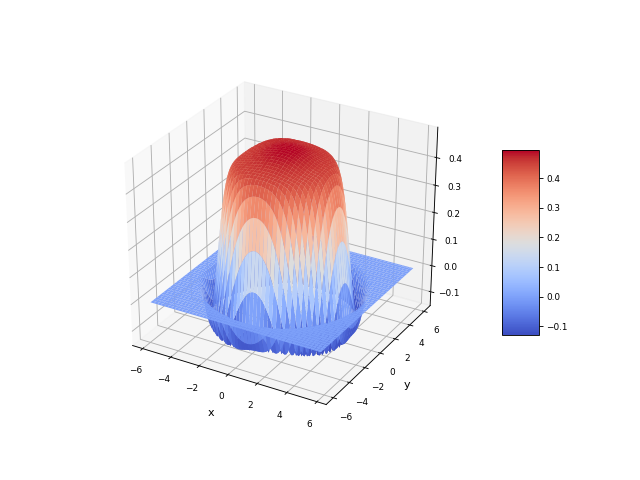}
               \caption{Learned $-\Delta q$}
    \end{subfigure}
    \hfill
    \begin{subfigure}[t]{0.47\textwidth}
             \centering
    \includegraphics[width = 1.2\textwidth]{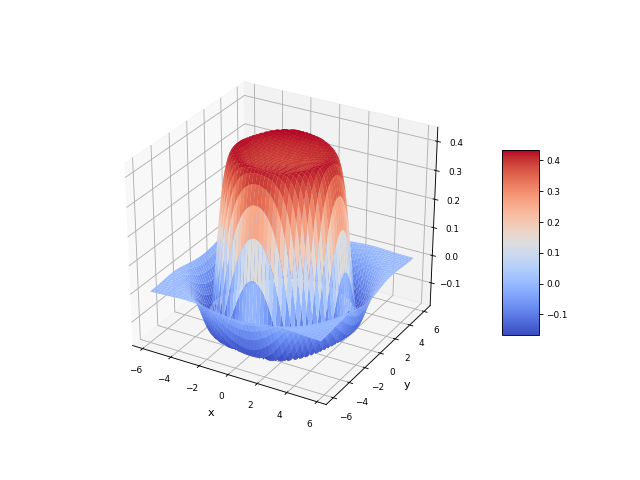}
               \caption{Learned $\displaystyle \partial_t \sigma$}
    \end{subfigure}
    \caption{\textbf{10D, $L^2-q-\sigma$  formulation \eqref{qsigma_full}:} predicted $q$, $\sigma$ and their partial derivatives.}
    \label{fig:qsigma_PDE_l2}
\end{figure}

\begin{figure}[htbp]
     \centering
     \begin{subfigure}[t]{0.3\textwidth}
         \centering
         \includegraphics[width=1.2\textwidth]{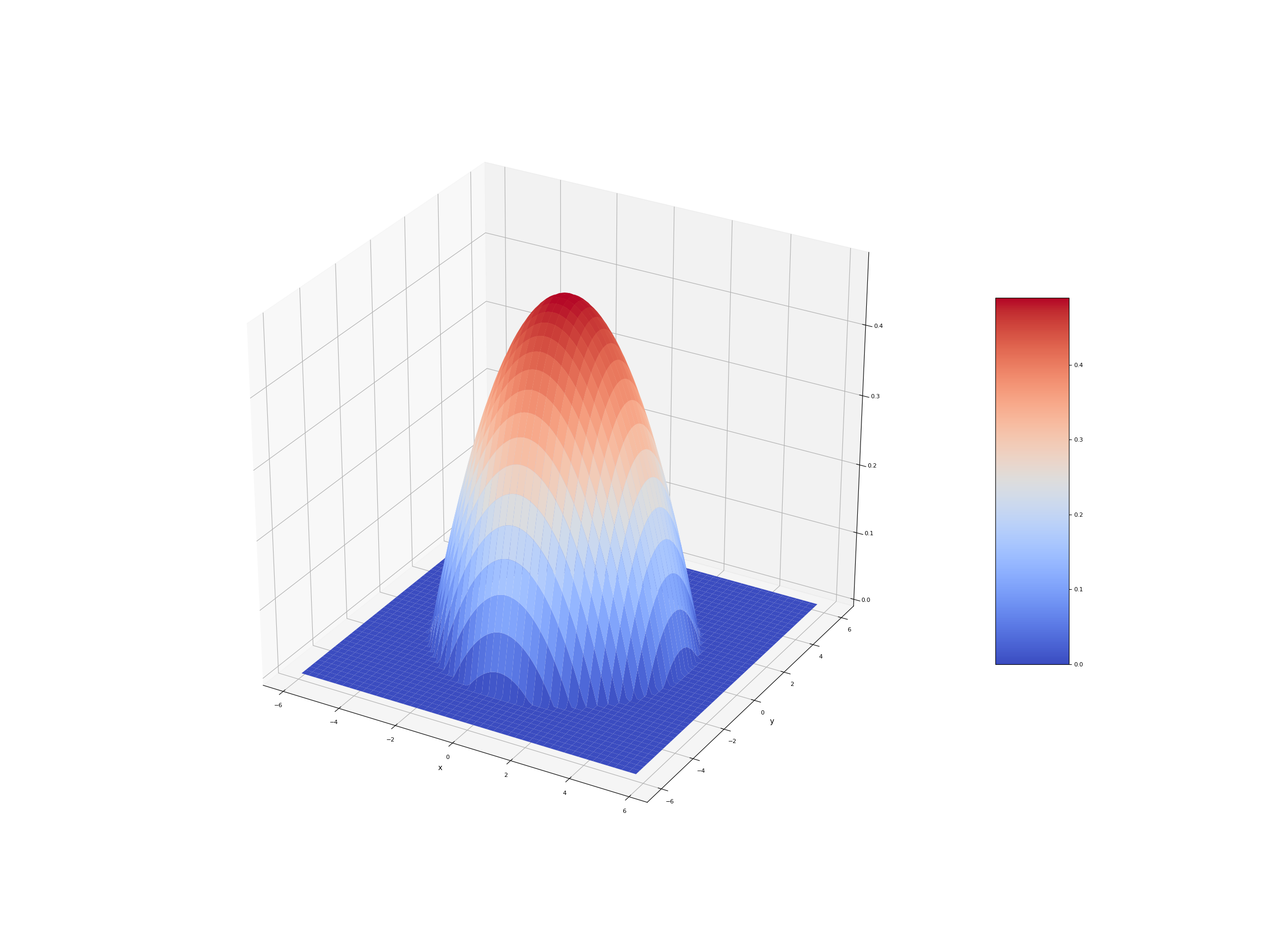}
         \caption{Barenblatt reference solution  }
     \end{subfigure}
    \hfill
     \begin{subfigure}[t]{0.3\textwidth}
         \centering
         \includegraphics[width=1.2\textwidth]{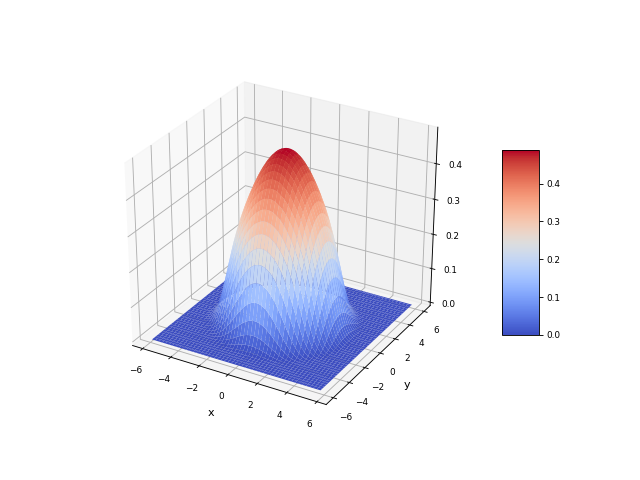}
         \caption{Learned solution slice}
     \end{subfigure}
     \hfill
    \begin{subfigure}[t]{0.3\textwidth}
         \centering
\includegraphics[width=1.2\textwidth]{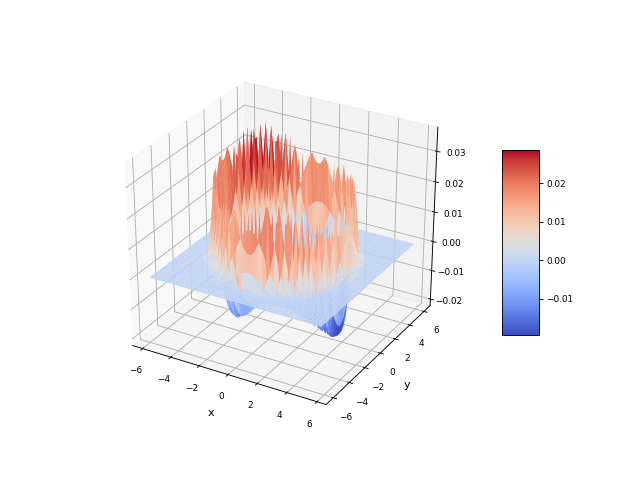}
         \caption{Learned solution error}
     \end{subfigure}
        \caption{\textbf{10D,  $L^1-q-\sigma$ formulation \eqref{qsigma_full}:} Predicted solution slice  $u(0.5,x,y,1.0,\cdots, 1.0)$ for $\x\in \T = [-6,6]^{10}$, $t= 0.5$. }
        \label{fig:qsigma_l1}
\end{figure}

\begin{figure}
    \centering
    \begin{subfigure}[t]{0.47\textwidth}
             \centering
    \includegraphics[width = 1.2\textwidth]{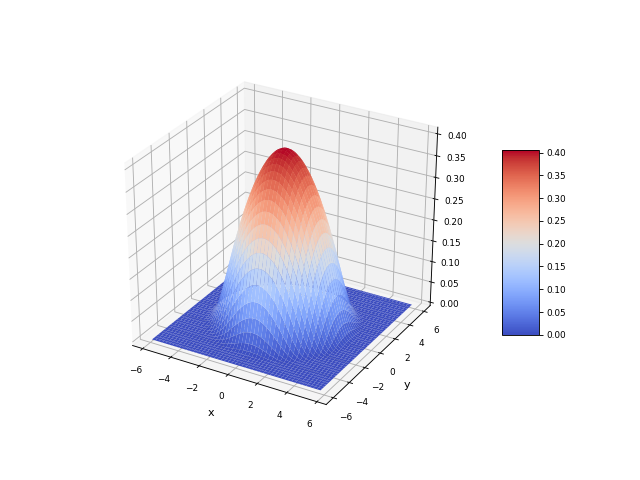}
               \caption{Learned $q$}
    \end{subfigure}
    \hfill
    \begin{subfigure}[t]{0.47\textwidth}
             \centering
    \includegraphics[width = 1.2\textwidth]{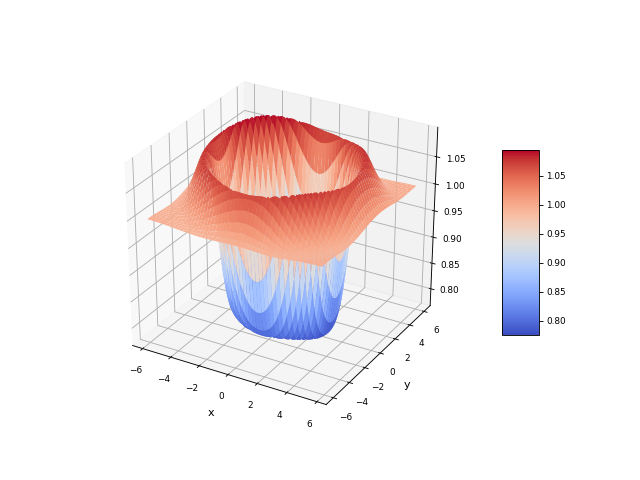}
               \caption{Learned $\displaystyle \sigma$}
    \end{subfigure}
    \begin{subfigure}[t]{0.47\textwidth}
             \centering
    \includegraphics[width = 1.2\textwidth]{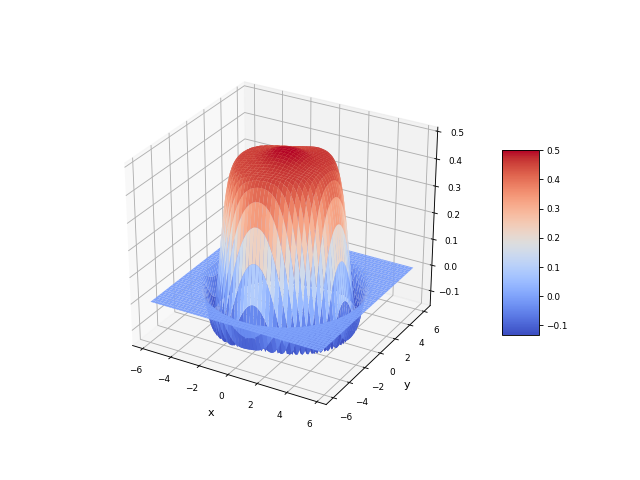}
               \caption{Learned $-\Delta q$}
    \end{subfigure}
    \hfill
    \begin{subfigure}[t]{0.47\textwidth}
             \centering
    \includegraphics[width = 1.2\textwidth]{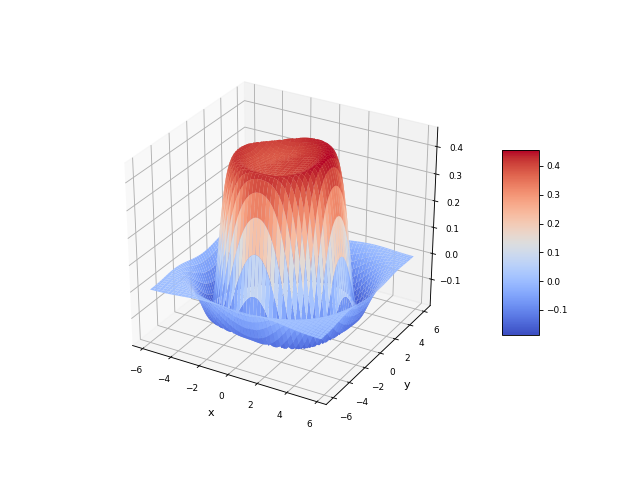}
               \caption{Learned $\displaystyle \partial_t \sigma$}
    \end{subfigure}
    \caption{\textbf{10D, $L^1-q-\sigma$  formulation \eqref{qsigma_full}:} predicted $q,\sigma$ and their partial derivatives.}
    \label{fig:qsigma_PDE_l1}
\end{figure}

\begin{figure}
    \centering
     \begin{subfigure}[t]{0.45\textwidth}
             \includegraphics[width =1.2\textwidth]{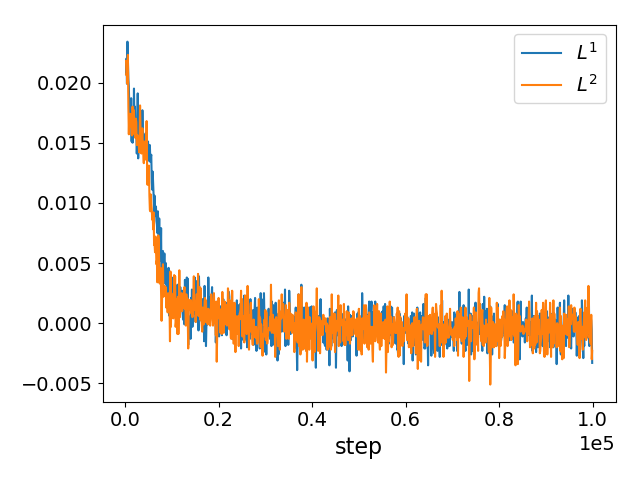}
             \caption{\textbf{5D}}
    \end{subfigure}
    \hfill
    \begin{subfigure}[t]{0.45\textwidth}
             \includegraphics[width =1.2\textwidth]{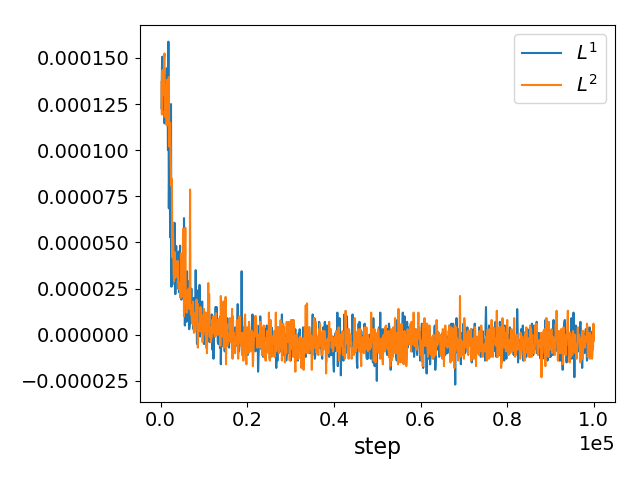}
             \caption{\textbf{10D}}
    \end{subfigure}
    \caption{Empirical $\int_Q U_2^2(t,\x) +\mathcal{L}_{q,\sigma}\big(u_{q,\sigma}(t,\x;\theta_{q},\theta_{\sigma} )\big)$ v.s. training steps. }
    \label{fig:qsigma_usq_comp}
\end{figure}

Specifically, the homogeneous Dirichlet boundary condition is imposed as a hard constraint following \eqref{q_with_conditions} and the condition for $\sigma$ is imposed with \eqref{sigma_condition}. 
The condition \eqref{sigma_selection} was not strongly imposed for training reasons. The initial condition can then be softly enforced with term $\mathcal{L}_{I}$ as mentioned earlier. The specific algorithmic settings are further presented in Table \ref{tab:qsigma_error} along with the relative errors computed for the trained solution slice $u_{q,\sigma}(0.5,x,y,1.0,\cdots,1.0; \theta_{q}, \theta_{\sigma})$ at time $t= 0.5$ comparing with the exact solution \eqref{exact}. The comparison of predicted solutions with exact solution are further presented in Figure \ref{fig:qsigma_l2} and Figure \ref{fig:qsigma_l1}. In addition, the predicted function $q$ and $\sigma$ are depicted in Figure \ref{fig:qsigma_PDE_l2} and Figure \ref{fig:qsigma_PDE_l1}. These figures are further used to show the predicted $-\Delta q$ and $\partial_t \sigma$ to verify that the condition
$$\Delta q + \partial_t \sigma = 0$$
is satisfied. Finally, Figure \ref{fig:qsigma_usq_comp} is used to demonstrate the computational value of $\mathcal{L}_{q,\sigma}$ converges to
to $-\int_{Q}U_2^2$ as training proceeds.
This observation confirms once again the theoretical result \eqref{form_equivalency} derived in \cite{brenier2020examples}.
Here, the batch of training data at each training step is used to evaluate empirically the value of  $-\int_Q U_2^2$ for the exact solution $U_2(t,\x)$ defined as in \eqref{exact}. Such value is further compared with the empirical loss $\mathcal{L}_{q,\sigma}(u_{q,\sigma})$ at that time. The difference of the these values gradually reduces as the training continues (Figure \ref{fig:qsigma_usq_comp}).

\begin{table}[htbp]
  \centering
    \begin{tabular}{|c|l|c|c|c|c|c|c|}
    \toprule
         \textbf{Dimension} &  & \multicolumn{1}{r|}{\textbf{1}} & \multicolumn{1}{r|}{\textbf{2}} & \multicolumn{1}{r|}{\textbf{3}} & \multicolumn{1}{r|}{\textbf{4}} & \multicolumn{1}{r|}{\textbf{5}} & \multicolumn{1}{r|}{\textbf{10}} \\
    \midrule
    \multirow{2}[4]{4cm}{\textbf{Relative Errors(\%) for $L^2-q-\sigma$ Formulation}} & \boldmath{}\textbf{$L^2$}\unboldmath{} & \multicolumn{1}{r|}{\textbf{1.95}} & \multicolumn{1}{r|}{\textbf{3.2}} & \multicolumn{1}{r|}{\textbf{3.88}} & \multicolumn{1}{r|}{\textbf{3.97}} & \multicolumn{1}{r|}{\textbf{4.77}} & \multicolumn{1}{r|}{\textbf{4.03}} \\
\cmidrule{2-8}          & \boldmath{}\textbf{$L^1$}\unboldmath{} & \multicolumn{1}{r|}{\textbf{1.64}} & \multicolumn{1}{r|}{\textbf{3.11}} & \multicolumn{1}{r|}{\textbf{3.5}} & \multicolumn{1}{r|}{\textbf{4.02}} & \multicolumn{1}{r|}{\textbf{5.11}} & \multicolumn{1}{r|}{\textbf{4.14}} \\
    \midrule
    \multirow{2}[4]{4cm}{\textbf{Relative Errors(\%) for $L^1-q-\sigma$ Formulation}} & \boldmath{}\textbf{$L^2$}\unboldmath{} & \multicolumn{1}{r|}{\textbf{2.06}} & \multicolumn{1}{r|}{\textbf{2.96}} & \multicolumn{1}{r|}{\textbf{3.83}} & \multicolumn{1}{r|}{\textbf{3.94}} & \multicolumn{1}{r|}{\textbf{5.63}} & \multicolumn{1}{r|}{\textbf{4.28}} \\
\cmidrule{2-8}          & \boldmath{}\textbf{$L^1$}\unboldmath{} & \multicolumn{1}{r|}{\textbf{1.72}} & \multicolumn{1}{r|}{\textbf{2.79}} & \multicolumn{1}{r|}{\textbf{3.28}} & \multicolumn{1}{r|}{\textbf{3.72}} & \multicolumn{1}{r|}{\textbf{6.41}} & \multicolumn{1}{r|}{\textbf{4.59}} \\
    \midrule
    \multirow{3}[6]{*}{Formulation Weights } & $\nu$ & \multicolumn{6}{c|}{1} \\
\cmidrule{2-8}          & $\kappa$ & \multicolumn{6}{c|}{$10^3$} \\
\cmidrule{2-8}          & $\gamma$ & \multicolumn{5}{c|}{$10^3$}             & \multicolumn{1}{r|}{1} \\
    \midrule
    \multirow{2}[4]{*}{NN Architecture} & Width/Depth & \multicolumn{6}{c|}{200/2} \\
\cmidrule{2-8}          & \# trainable & \multicolumn{1}{r|}{82002} & \multicolumn{1}{r|}{82402} & \multicolumn{1}{r|}{82802} & \multicolumn{1}{r|}{83202} & \multicolumn{1}{r|}{83602} & \multicolumn{1}{r|}{85602} \\
    \midrule
    \multirow{2}[4]{*}{Data Sampling} & $\theta_0$ & \multicolumn{6}{c|}{0.3} \\
\cmidrule{2-8}          & $\theta_1$ & \multicolumn{6}{c|}{0.3} \\
    \midrule
    \multirow{2}[4]{*}{Training} &   Steps & \multicolumn{6}{c|}{$10^5$} \\
\cmidrule{2-8}          & Learning Rate & \multicolumn{6}{c|}{$10^{-3}$} \\
    \bottomrule
    \end{tabular}%
    \caption{\textbf{$q-\sigma$ formulation \eqref{qsigma_form} (hard Dirichlet B.C.+soft I.C.): } Relative error comparison for various dimensions.}
  \label{tab:qsigma_error}%
\end{table}%

\section{Numerical example: waiting-time phenomena }\label{sec:nuemrical_waiting}
In this section, we further consider the following IBVP
\begin{equation}\label{numerical_waiting}
\begin{split}
     &\partial_t u =  \frac{1}{2}\Delta u^2\quad (t,\x)\in Q = [0,1]\times\T,\\
     &u(0,\x) = u_0(\x) =\begin{cases} 
      \cos(|\x|), &  |\x|\leq \frac{\pi}{2} \\
      0, & \text{elsewhere}
      \end{cases}\\
      &u(t,\x)|_{\partial \T} =0.
\end{split}
\end{equation}
where $\T = [-4,4]^d$.
The general exact solution to \eqref{numerical_waiting} can hardly be derived. When $d=2$, the reference solution can be taken to be the numerical solutions obtained with a moving mesh finite element method following \cite{wells2004moving} instead. In particular, the mesh is advanced forward in time using a forward Euler time-stepping scheme. These mesh-based results are then compared with the ones obtained following a deep learning framework under various formulations.
For higher dimensions, the mesh-based solver in general will suffer from curse of dimensionality. The moving mesh method would also be more challenging to design. Therefore, for comparison reasons, we only present the results for $d=2$ while noticing the higher dimensional cases can also be handled by the neural network based algorithms.

\begin{figure}[htbp]
    \centering
     \begin{subfigure}[b]{\textwidth}
    \includegraphics[width = \textwidth]{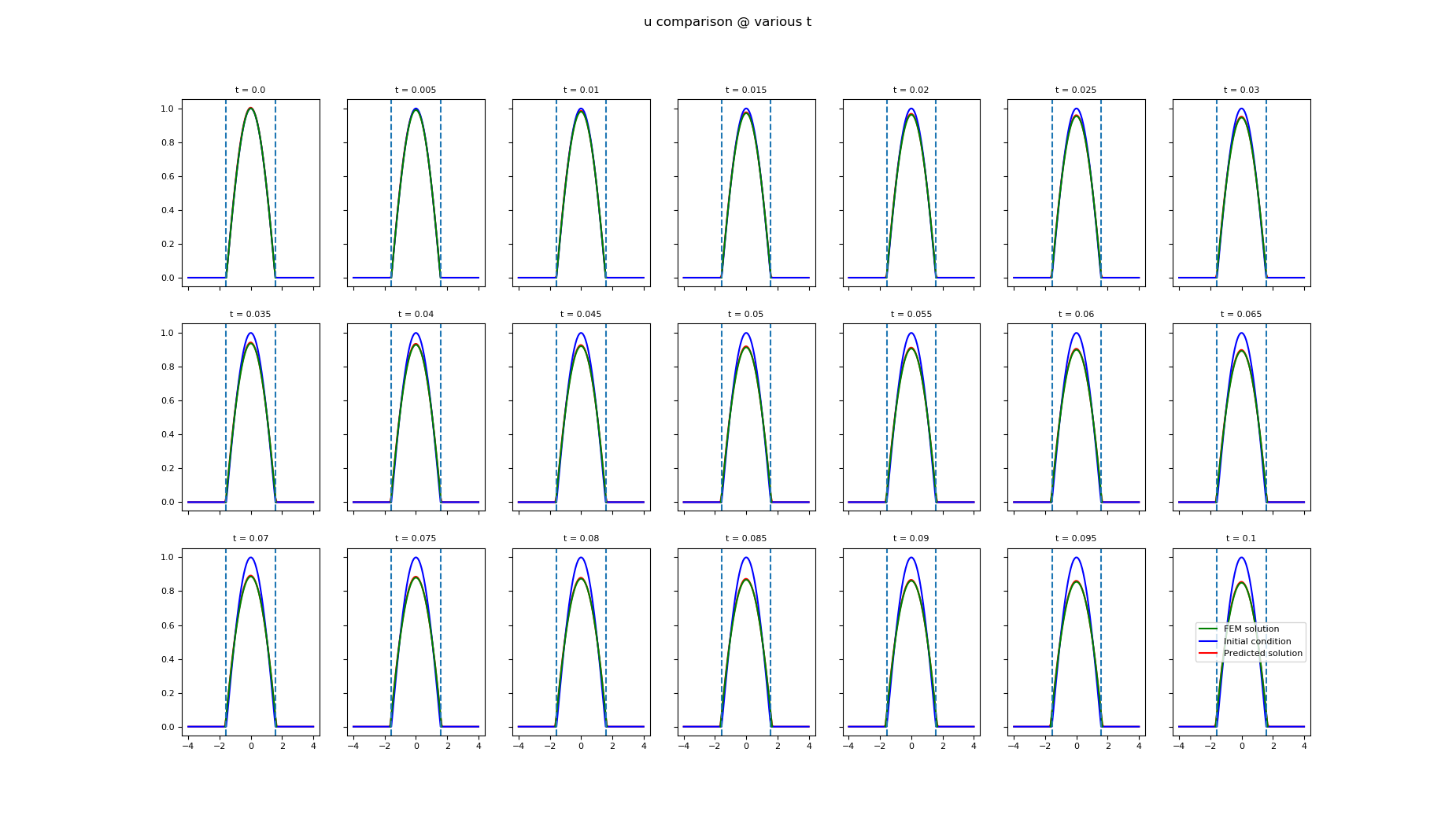}
    \caption{Snapshot solution slices for $t\in [0,0.1]$}
    \label{waiting:PINN_a}
    \end{subfigure}
     \begin{subfigure}[b]{\textwidth}
    \includegraphics[width = \textwidth]{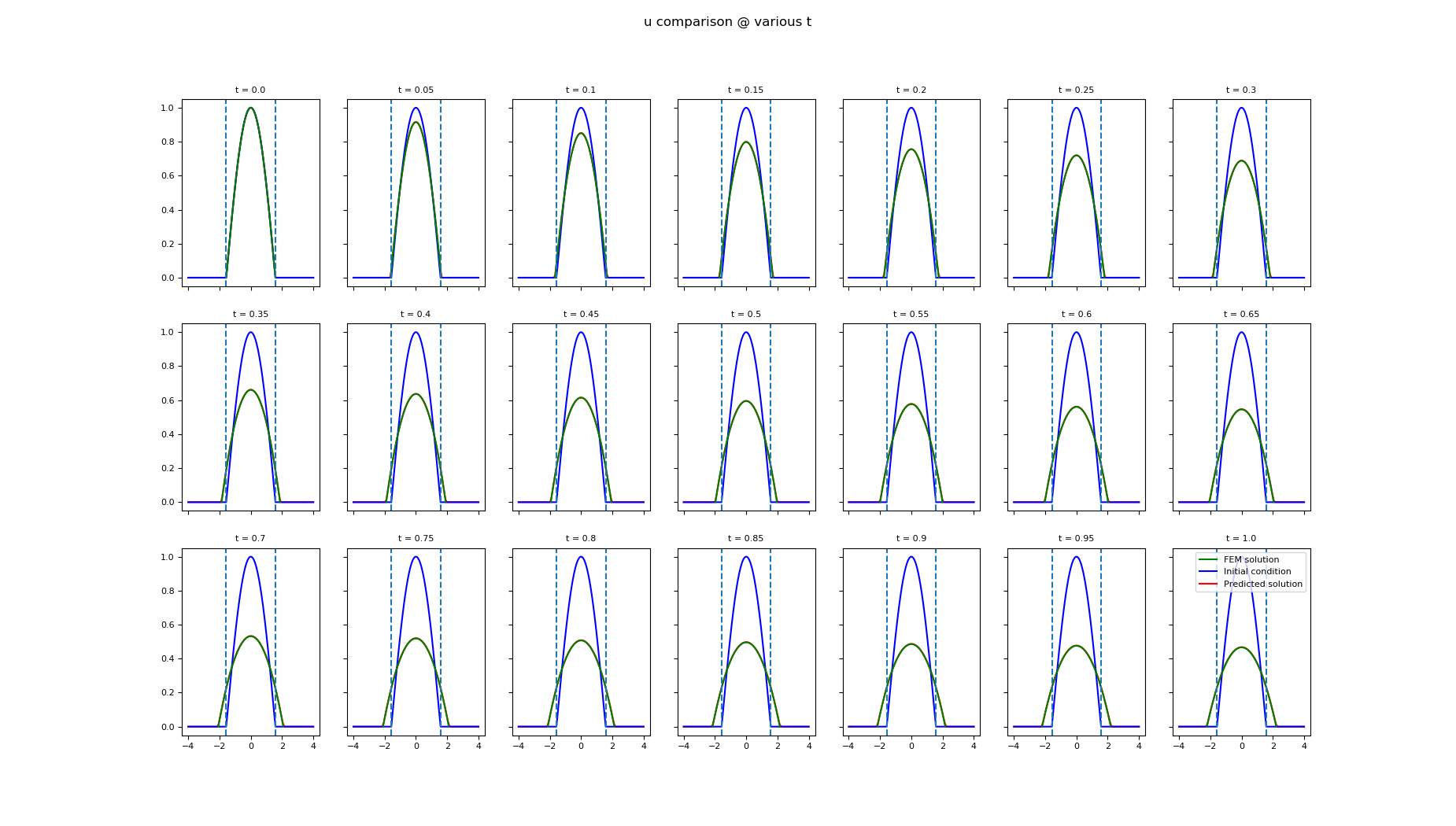}
    \caption{Snapshot solution slices for $t\in [0,1.0]$}
        \label{waiting:PINN_b}
    \end{subfigure}
    \caption{\textbf{2d, $L^2-$ PINN formulation, waiting-time phenomena \eqref{numerical_waiting}:} snapshots of solution cross-section $u(t,x,0)$ at $y=0$. Green: reference solutions obtained with moving mesh FEM ($DOF= 901$); red: predicted solutions obtained with PINN formulation; blue: initial condition.}
    \label{fig: PINN_waiting}
\end{figure}

\begin{figure}[htbp]
    \centering
     \begin{subfigure}[b]{\textwidth}
    \includegraphics[width = \textwidth]{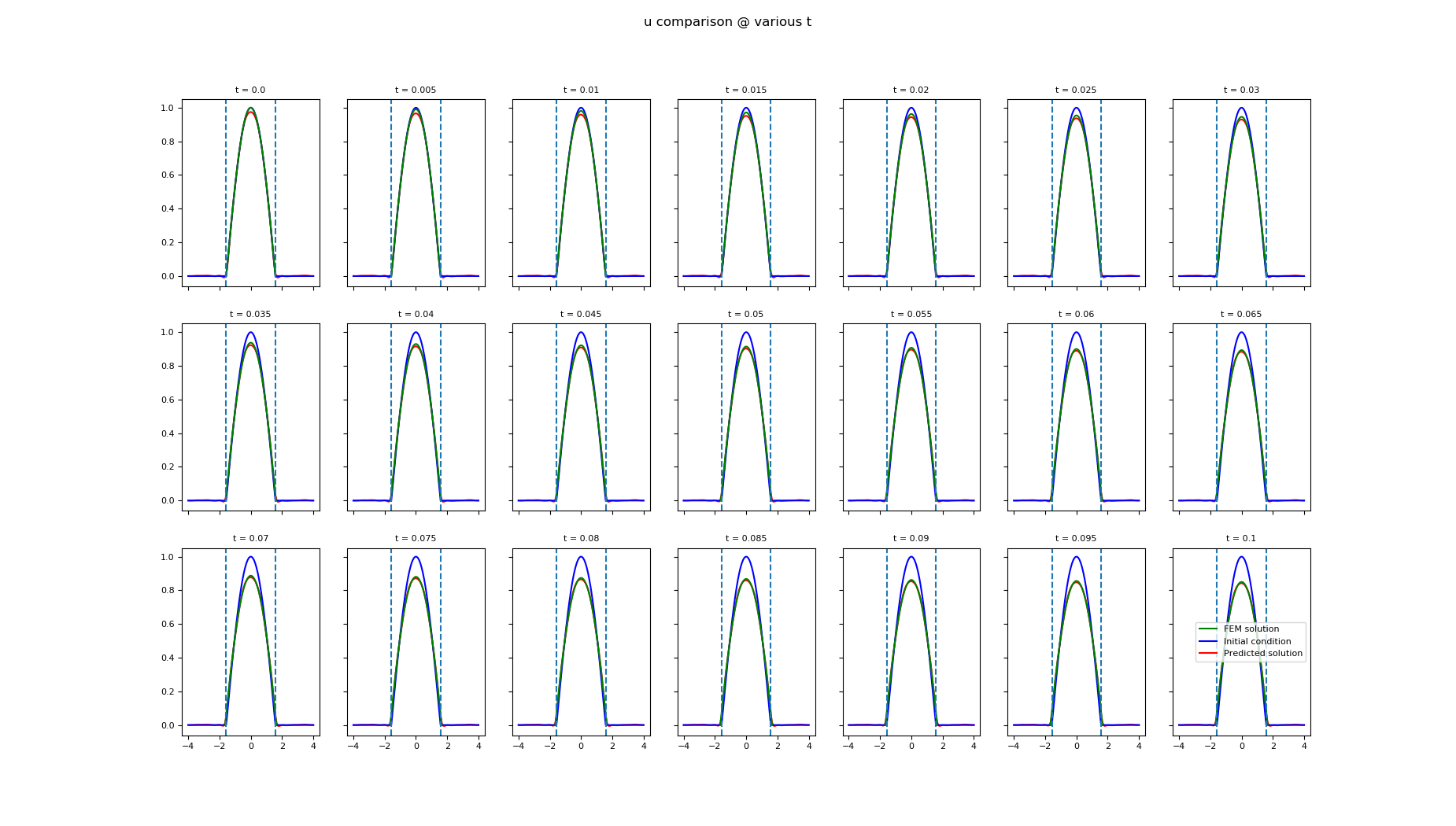}
    \caption{$t\in [0,0.1]$}
    \label{waiting:phi_a}
    \end{subfigure}
     \begin{subfigure}[b]{\textwidth}
    \includegraphics[width = \textwidth]{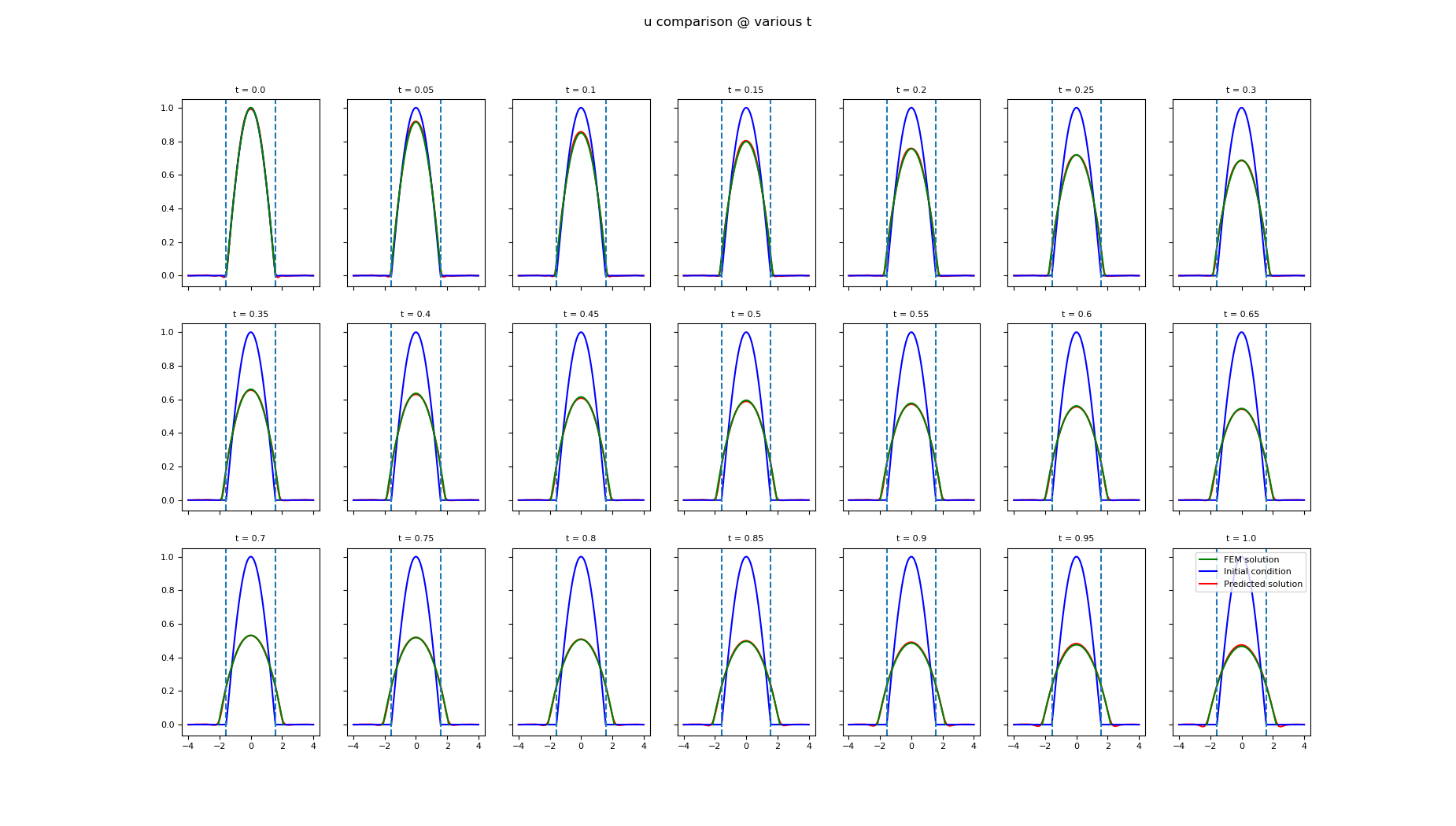}
    \caption{Snapshot solution slices for $t\in [0,1.0]$}
        \label{waiting:phi_b}
    \end{subfigure}
    \caption{\textbf{2d, $L^2-\phi$  formulation, waiting-time phenomena \eqref{numerical_waiting}:} snapshots of solution cross-section $u(t,x,0)$ at $y=0$. Green: reference solutions obtained with moving mesh FEM ($DOF= 901$); red: predicted solutions obtained with $\phi$ formulation; blue: initial condition.}
    \label{fig: phi_waiting}
\end{figure}

We also note that the solution to PME of this type of initial condition exhibits a waiting-time phenomenon \cite{ngo2017study}. In fact, the velocity of the free boundary of QPME is given by Darcy's law \cite{shmarev2005interfaces}, i.e.
\begin{equation}\label{darcy}
    \Gamma'(t) = \lim_{\x\to \Gamma(t)^-}\nabla (\frac{u^2}{2}),
\end{equation}
where the limit is taken from the interior of the support.
Thus, as one may compute, the free boundary of solution to 
\eqref{numerical_waiting} should not move until a finite amount of time has elapsed as initially $\Gamma'(0) = \nabla (\frac{u^2_0}{2})$ vanishes at the free boundary $\Gamma(0): |\x| = \frac{\pi}{2}$. This phenomena of waiting can be observed form solutions obtained with Finite Element Method as shown in Figure \ref{fig: PINN_waiting}, where the dashed vertical lines indicate the exact initial location of the free boundary. In Figure \ref{waiting:PINN_a}, a series of snapshots of the solution for $t\in[0,0.1]$ is plotted, while in \ref{waiting:PINN_b}, the solution snapshots for a broader range of time is plotted. As one may observe, the free boundary of the solution barely moves in the entire time of $t\in [0,0.1]$ and only start to change by the time of $t= 0.2$.  This phenomena can also be accurately captured by a solution obtained following the PINN formulation \eqref{PINN_full}. In specific, the solution obtained with the PINN formulation is presented in comparison with the ones obtained with the moving mesh FEM in Figure \ref{fig: PINN_waiting}. The solution slices essentially overlap one another. Parallelly, problem \eqref{numerical_waiting} is also solved with the $\phi$ formulation \eqref{phi_form}. The comparison  of the resulted solution with a FEM solution is then presented as in Figure \eqref{fig: phi_waiting}, which verifies the effectiveness of the neural network based solutions using $\phi$ formulation. In these numerical tests, the choices of the hyper-parameters are taken to be the same as in Table \ref{tab:PINN_error} and Table \ref{tab:phi_error} respectively. 

 \section{Conclusion} 
 In this paper, we explored different variational forms in solving high-dimensional QPME with neural networks. In specific, three formulations were considered. For the PINN formulation, the solution is directly parametrized and the PDE residual is minimized. A theoretical analysis is carried out to show that the convergence of the PINN loss guarantees a convergence to the exact solution in $L^1$ sense. Moreover, this analysis also suggests the use of the $L^1$ norm to quantify the residual and approximation error. 
 In addition, inspired by the work  \cite{brenier2020examples},  a $\phi$ formulation and a $q-\sigma$ formulation is further presented and used to solve the QPME in a very weak sense. Theoretically, these formulation can identify solutions with less regularity. All formulations are then tested with the Barrenbaltt solution in low and high dimensions. Experiments have shown that $\phi$ formulation and $q,\sigma$ formulation can provide approximate solutions with a similar level of accuracy compared with PINN in low-dimensional cases but the optimization aspect continues to pose challenges in high-dimensional cases. A two-dimensional example of QPME that exhibits waiting phenomena is also presented to show the capability of deep learning based methods in identifying solution features as such.  
 
 Other aspects of the discussion toward solving QPME with deep learning includes the hard and soft imposition of certain conditions of the solutions in all formulations such as initial conditions and boundary conditions. Additionally, an efficient sampling scheme is proposed aiming at a faster convergence towards the solution desired especially in high-dimensional cases. These treatments in principal can also be applied in other scenarios where the PDE solutions are parametrize with neural networks.

While such efforts can all contributes to more efficient implementation of solving high-dimensional QPMEs, we must admit that the training success is overwhelmed by the large number of hyper-parameters. Moreover, for practical applications,
neural network training using stochastic gradient descent type schemes, which means one must accept a significant and unavoidable uncertainty about optimization success.
An efficient strategy on making choices of hyper-parameters could potentially be an interesting direction for future investigations.
More broadly, whether a similar variational form could be derived for general $m$ of porous medium equation is also an open question.

\bibliographystyle{unsrt}
\bibliography{reference}
\end{document}